\documentclass[11pt]{article}
\usepackage{IEEEtrantools}
\usepackage{mathtools, nccmath}
\usepackage{latexsym}
\usepackage{amsfonts}
\usepackage{amssymb}
\usepackage{psfrag}
\usepackage{graphicx}
\usepackage{epsfig}
\usepackage{amsmath,amsfonts,amsthm}
\usepackage{mathrsfs}
\usepackage{enumerate}
\usepackage{enumitem} 
\usepackage{dsfont}
\usepackage{hyperref}
\usepackage{pst-node}
\usepackage{cleveref}
\usepackage{enumitem}
\usepackage{tikz-cd}
\usepackage{multirow}
\usepackage{tikz}
\usetikzlibrary{arrows,decorations.markings}
\usepackage{bookmark}
\usepackage{hyperref}
\hypersetup{
     colorlinks   = true,
     citecolor    = blue
}

\newcommand\restr[2]{{
  \left.\kern-\nulldelimiterspace 
  #1 
  \littletaller 
  \right|_{#2} 
  }}

\newcommand{\littletaller}{\mathchoice{\vphantom{\big|}}{}{}{}}

\newcommand{\fn}{\mathcal{F}_n}

\newcommand{\be}{\begin{equation}}
\newcommand{\ee}{\end{equation}}
\newcommand{\bea}{\begin{eqnarray}}
\newcommand{\eea}{\end{eqnarray}}
\newcommand{\bean}{\begin{eqnarray*}}
\newcommand{\eean}{\end{eqnarray*}}
\newcommand{\brray}{\begin{array}}
\newcommand{\erray}{\end{array}}
\newcommand{\biearray}{\begin{IEEEarray}{rCl}}
\newcommand{\eiearray}{\end{IEEEarray}}

\newtheorem{dfn}{Definition}[section]
\newtheorem{thm}[dfn]{Theorem}
\newtheorem{lmma}[dfn]{Lemma}
\newtheorem{ppsn}[dfn]{Proposition}
\newtheorem{crlre}[dfn]{Corollary}
\newtheorem{xmpl}[dfn]{Example}
\newtheorem{rmrk}[dfn]{Remark}

\newcommand{\bdfn}{\begin{dfn}\rm}
\newcommand{\bthm}{\begin{thm}}
\newcommand{\blmma}{\begin{lmma}}
\newcommand{\bppsn}{\begin{ppsn}}
\newcommand{\bcrlre}{\begin{crlre}}
\newcommand{\bxmpl}{\begin{xmpl}}
\newcommand{\brmrk}{\begin{rmrk}\rm}

\newcommand{\edfn}{\end{dfn}}
\newcommand{\ethm}{\end{thm}}
\newcommand{\elmma}{\end{lmma}}
\newcommand{\eppsn}{\end{ppsn}}
\newcommand{\ecrlre}{\end{crlre}}
\newcommand{\exmpl}{\end{xmpl}}
\newcommand{\ermrk}{\end{rmrk}}

\newcommand{\bbc}{\mathbb{C}}

\newcommand{\bbn}{\mathbb{N}}

\newcommand{\tr}{\mathrm{tr}}

\def \qed { \mbox{}\hfill
$\Box$\vspace{1ex}}

\title{Sections and Chapters}

\begin{document}
	
	
	\author{\sc{Keshab Chandra Bakshi, Satyajit Guin, Biplab Pal, Sruthymurali}}
	\title{Higher reflections and entropy of canonical shifts for inclusions of $C^*$-algebras with finite Watatani index}	\maketitle
	
	
\begin{abstract}
Given a unital inclusion of simple $C^*$-algebras equipped with a conditional expectation of index-finite type, we study  Fourier transforms and rotation operators and introduce the reflection operators on the relative commutants. We prove that the reflections are unital, involutive, $*$-preserving anti-homomorphisms that preserve certain Markov-type traces. As an application, we prove Fourier theoretic inequalities on the higher relative commutants and refine the existing constant in Young’s inequality as presented in the current literature. By employing the reflection operators, we define a canonical shift on the von Neumann algebra generated by the relative commutants. We establish a connection between the Connes–Størmer entropy of the canonical shift and the minimal Watatani index.
\end{abstract}
\bigskip
	
	{\bf AMS Subject Classification No.:} {\large 46}L{\large 37}\,, {\large 47}L{\large 40}\,, {\large 46}L{\large 05}\,, {\large 43}A{\large 30}\,.
	
	\bigskip
	\hypersetup{linkcolor=blue}
	\tableofcontents
	\bigskip

\section{Introduction}
Vaughan Jones in 1983 (see \cite{Jo}) invented a notion of index $[M:N]$ for a subfactor $N\subset M$ of type $II_1$ factors and proved that index, rather surprisingly,  takes value in the set $\{4\cos^2 \pi/n: n>3\}\cup [4,\infty]$. In another direction, in a remarkable paper (see \cite{CS}) Connes and St{\o}rmer had discovered a notion of  entropy of the $*$-automorphism (later modified by Choda for $*$-endomorphism in \cite{MCJOT}) of a finite von Neumann algebra as a noncommutative generalization of Kolmogorov-Sinai entropy. In that same paper, Connes and St{\o}rmer also introduced a notion of noncommutative relative entropy $H(B_1|B_2)$ for a pair of finite dimensional $C^*$-subalgebras $B_1$ and $B_2$ of a finite von Neumann algebra $M$ as a noncommutative generalization of  the classical conditional entropy in ergodic theory. Pimsner and Popa in \cite{PP} observed that the definition of the relative entropy works for arbitrary von Neumann subalgebras $B_1, B_2\subset M$. They had discovered a surprising connection between the Connes-St{\o}rmer relative entropy $H(M|N)$ and the Jones index $[M:N]$, namely, a subfactor $N\subset M$ is {extremal} if and only if $H(M|N)=\log [M:N].$ In particular, if a subfactor $N\subset M$ has finite depth, then $H(M|N)=\log [M:N]$. They have also provided a forumula for $H(M|N)$ in the general situation in terms of the atoms of the relative commutant $N^{\prime}\cap M$ which is a finite dimensional $C^*$-algebra consisting of the elements of $M$ that commutes with $N$. 

In the landmark paper \cite{Jo}, Jones also  introduced the so-called {basic construction tower} of $II_1$ factors $N\subset M\subset M_1\subset \cdots \subset M_k \subset  \cdots $ that plays a significant role in the index theory of subfactors. The finite dimensional $C^*$-algebras $N^{\prime}\cap M_k$ and $M^{\prime}\cap M_l$ are called the relative commutants of the subfactor $N \subset M$, which play a fundamental role in subfactor theory. The relative commutants were crucially used by Jones in \cite{Jo2} to develop his celebrated theory of  planar algebras which are complete invariant for a `\textit{good class of subfactors'}. In formalism of the structure of Jones' planar algebras, the most interesting algebraic ingredients were the Fourier transforms and rotation maps on the relative commutants. These maps were also used by Bisch to investigate intermediate subfactors in \cite{Bi94}. Recently, as an important application of the the Fourier transform and rotation maps, Jiang, Liu and Wu \cite{Liunoncommutative} provided a non-commutative version of the Hausdorff-Young inequality, the Young’s inequality, and various uncertainty principles on the second relative commutants $N^{\prime}\cap M_1$ and $M^{\prime}\cap M_2$ of a subfactor $N\subset M$ with $[M:N]<\infty$ using planar algebraic tools. In another direction, to classify subfactors of the hyperfinite $II_1$ factors with small Jones index, Ocneanu had introduced an important endomorphism on the tower of relative commutants called the `\textit{canonical shift}' on the relative commutants. The canonical shift, usually denoted by  $\Gamma$, naturally induces a 2-shift on the injective von Neumann algebra generated by the tower of the relative commutants. The canonical shift has strong resemblance with Longo's canonical endomorphism \cite{Lo2,Lo1}. Motivated by Popa’s remarkably deep result, where he established a criterion for strong amenability of subfactors in terms of relative entropy \cite{PAmenable}, Choda \cite{Choda} discovered that the Connes-St{\o}rmer  entropy of the canonical shift is closely related to the Jones index. More explicitly, she proved that if a subfactor has generating property and is extremal, then $H(\Gamma)= \log [M:N].$ In particular, for a  subfactor with  finite depth, we have
$$H(\Gamma)= H(M|N)= \log [M:N].$$ Further investigation in this direction has been done by Hiai (\cite{Hia,Hia2}). These results on canonical shifts and entropy for type $II_1$ factors have also been extended to the type $III$ factors \cite{CH}. Recall that, analogous to  the Jones index, Kosaki \cite{Ko} introduced a notion of index and basic construction for infinite factors. The theory of type $III$ factors have been studied further by Popa, Longo, Izumi, Loi and many others (for instance, see \cite{Izu2,Loi,Lo3,Lo4,PAmenable}). In \cite{CH}, Choda and Hiai introduced the canonical shift for a type $III$-subfactor with a conditional expectation. They found the relationship between the Kosaki index and entropy, like in the case of $II_1$ subfactors.

Inclusion of simple $C^*$-algebras encompasses both type  $II_1$ and type $III$-subfactors. Due to the recent advances in the $C^*$-algebra classification program, the study of the inclusions of $C^*$-algebras has become an active area of research (see \cite{BakshiVedlattice,SieRor,IzumiNew,Mukohara,Ror}, to name a few). Generalizing both the Jones index (for type $II_1$-subfactors) and Kosaki index (for type $III$-subfactors), Watatani proposed an algebraic notion of index of a conditional expectation having a `\textit{quasi-basis}' \cite{Watataniindex} for an inclusion of unital $C^*$-algebras $B\subset A$. For  an inclusion of simple unital $C^*$-algebras $B\subset A$ with a conditional expectation of index-finite type, he further proved that there is a `\textit{minimal conditional expectation'} $E_0$ and the \textit{minimal index} for the inclusion is defined as the Watatani index of $E_0$, denoted by ${[A:B]}_0$. The tower of simple $C^*$-algebras $B\subset A\subset A_1\subset \cdots \subset A_n\subset \cdots$ is the iteration of Watatani's $C^*$-basic construction (with respect to the {minimal conditional expectation}) which is analogous to the Jones tower of basic construction. The relative commutants $B^{\prime}\cap A_k$ and $A^{\prime}\cap A_l$ come equipped with a consistent Markov-type trace as explained in \cite{BakshiVedlattice}. Motivated by subfactor theory, the first author and Gupta introduced (in \cite{BakshiVedlattice}) the Fourier transform $\mathcal{F}:B^{\prime}\cap A_1\to A^{\prime}\cap A_2$, the rotation map ${\rho}_{+}:B^{\prime}\cap A_1\to B^{\prime}\cap A_1$ and the convolution product $*$ on $B^{\prime}\cap A_1$. By a well-known crossed product construction, this generalizes the Fourier theory for finite groups. More generally, using the Fourier transform $\mathcal{F}_n: B^{\prime}\cap A_{n} \rightarrow A^{\prime}\cap A_{n+1}$, we can conclude that  for each $n$, $B^{\prime}\cap A_n$ and $A^{\prime}\cap A_{n+1}$ are isomorphic as vector spaces. The $C^*$-Fourier theory has found significant applications in various contexts. For instance, in \cite{BDLRTR} and \cite{BakshiVedlattice}  the Fourier theory is applied to introduce a notion of angle between intermediate subalgebras of an inclusion of simple $C^*$-algebras and certain rigidity of the angle has been exploited to answer a question of Longo \cite{LONCMP} regarding the upper bound of the cardinality of the lattice of intermediate subalgebras (see also \cite{BGJBound}). Employing the Fourier theory in \cite{BGP}, a calculable formula is obtained for the von Neumann entropy of the (Fourier) dual angle operator associated with a broad class of quadruples of simple $C^*$-algebras. 

In this paper, we give a comprehensive study of the higher rotation operators $\rho^{\pm}_n$. Motivated by planar algebra theory, we call $B^{\prime}\cap A_{n-1}$ (resp. $A^{\prime}\cap A_n$) the positive $n$-box space $\mathcal{P}_{n,+}$  (resp. negative $n$-box space $\mathcal{P}_{n,-}$).  It was proved in \cite{BakshiVedlattice,BGS2022} that the rotation operator $\rho_1^{\pm}$ on $B^{\prime}\cap A_1$ (resp., $A^{\prime}\cap A_2$) is a trace preserving involutive anti-homomorphism. However, this is no longer true for higher rotations (\Cref{see}). Indeed, in \Cref{high rho+}, we have proved that ${(\rho_n^{\pm})}^{n+1}=id_{\mathcal{P}_{n+1,\pm}}$ and therefore, they are no longer involutive. Furthermore, one can easily check that the rotation operators are neither $*$-preserving nor anti-homomorphism. In this paper, we introduce reflection operators (that is, involutive) as follows. For $n\geq 0$, define the reflection operators $r^+_{2n+1}: B^\prime \cap A_{2n+1} \to B^\prime \cap A_{2n+1}$ by $r^+_{2n+1}=(\rho^+_{2n+1})^{n+1}$, and the dual reflection operator  $r^-_{2n+1}: A^\prime \cap A_{2n+2} \to A^\prime \cap A_{2n+2}$ by $r^-_{2n+1}=(\rho^-_{2n+1})^{n+1}.$ These opeartors are particularly well-behaved.
\smallskip

\noindent\textbf{Thorem A:} (See  Theorems \ref{antihomo}, \ref{r+ trace preserving}, and \ref{r- every property}) Let $B\subset A$ be an inclusion of simple unital $C^*$-algebras with a conditional expectation of index-finite type. Then, for each $n\geq 0$, the following results hold:
\begin{enumerate}
\item[$(i)$] The reflection operator $r^{+}_{2n+1}: B^{\prime }\cap A_{2n+1} \to B^\prime \cap A_{2n+1}$, is a unital, involutive, $*$-preserving anti-homomorphism. It also preserves the trace.
\item[$(ii)$] The reflection operator $r^-_{2n+1}: A^\prime \cap A_{2n+2} \to A^\prime \cap A_{2n+2}$, is a unital, involutive, $*$-preserving anti-homomorphism. It also preserves the trace.
\end{enumerate} 
\smallskip

Using the techniques developed in \cite{BakshiVedlattice}, Bakshi, Guin and Sruthymurali in \cite{BGS2022} proved a noncommutative version of the Hausdorff-Young inequality, and a couple of uncertainty principles in the setting of inclusion of simple $C^*$-algebras for the positive $2$-box space $\mathcal{P}_{2,+}$, and thereby generalized a few results of \cite{Liunoncommutative}. In this paper, using higher Fourier transform we generalize these to higher $n$-box space, both positive and negative.  As an important application of Theorem A, we improve the existing constant in the Young's inequality obtained for two-box space in \cite{BGS2022}, and also prove the Young's inequality for negative two box space.
\smallskip

\noindent\textbf{Theorem B:} (See \Cref{Young}) Suppose $B\subset A $ is an inclusion of simple unital $C^*$-algebras with a conditional expectation of index-finite type. \textcolor{black}{Then, we have the following\,:}
\begin{enumerate}
     \item[$(i)$] For any $x,y \in B' \cap A_1$, we have
     $\,\lVert x * y\rVert _r \leq\dfrac{\delta}{\kappa^{+}_{0}}\lVert x \rVert _p \lVert  y\rVert _q$,
     where $1\leq p,q,r\leq \infty$ and $\frac{1}{p}+\frac{1}{q}= \frac{1}{r}+1$.
     \item[$(ii)$] For any $w,z \in A' \cap A_2$, we have
     $\,\lVert w * z\rVert _r \leq\dfrac{\delta}{\kappa^{-}_{0}}\lVert w \rVert _p \lVert  z\rVert _q$,
     where $1\leq p,q,r\leq \infty$ and $\frac{1}{p}+\frac{1}{q}= \frac{1}{r}+1$.
     \end{enumerate}
\smallskip

As a  further application of reflection operators, we introduce an analogue of the canonical shift $\Gamma$ on the injective von Neumann algebra $P$ generated by the relative commutants of the inclusion $B\subset A$ of simple unital $C^*$-algebras with a conditional expectation of index-finite type. We show that this also induces a $2$-shift on $P$. We have obtained a relationship betwen the Connes-St{\o}rmer entropy of $\Gamma$ and Watatani's minimal index for a finite depth inclusion $B\subset A$ and this generalizes the subfactor case.

\noindent\textbf{Theorem C:} (See \Cref{entropy 2})
 Let  $B \subset A $ be an inclusion of simple unital $C^*$-algebras with a conditional expectation of index-finite type and assume that the inclusion is of finite depth. Then, 
$$\frac{1}{2}H_{\tr}(P| \Gamma(P))= H_{\tr}(\Gamma)= \log[A:B]_0 .$$

The subtle difference between results in this paper and the corresponding results in the subfactor theory lies in the fact that, unlike for finite factors, we neither have a tracial state on the $C^*$-algebra to begin
with nor the ‘modular conjugation operator’. We do not have nice pictorial descriptions of the Fourier transforms and rotations maps also. However, we have found a way around using the relationship between quasi-basis and minimal conditional expectations. Outline of the paper is as follows. In \Cref{preliminaries}, we review the necessary background, including Watatani’s algebraic framework for $C^*$-index theory and the concept of noncommutative entropy. \Cref{Higher rotation operators} focuses on the study of higher rotation and reflection operators, where we derive several key formulae and investigate fundamental properties of the canonical shift. In \Cref{Fourier transforms revisited}, we develop a series of Fourier-theoretic inequalities for higher $n$-box spaces, with particular emphasis on an improved version of Young’s inequality. Finally, in \Cref{final section}, we establish a relationship between the Connes–Størmer entropy of the canonical shift and the minimal Watatani index.
 
\section{Preliminaries}\label{preliminaries}

In this section, we review some key concepts used throughout the paper. We avoid delving too deeply into these concepts; however, appropriate references are provided in places for the readers interested in digging deeper.

\subsection{$C^*$-index theory}
Let $B\subset A$ be an inclusion of simple unital $C^*$-algebras with a conditional expectation of index-finite type. Let the tower of $C^*$-basic construction for the inclusion $B\subset A$ be:
	\begin{center}
		$B \subset A \subset A_1 \subset A_2 \subset \cdots \subset A_n \subset \cdots$
	\end{center}
	equipped with the unique (dual) minimal conditional expectations $E_{n}: A_{n}\to A_{n-1}$, $n \geq 0$, where we adopt the conventions $A_{-1}:=B$ and $A_0:=A$. For each $n\geq 1$, let $e_n $ be the Jones projection in $A_n$ implementing the $C^*$-basic construction of $A_{n-2} \subset A_{n-1}$ corresponding to the minimal conditional expectation $E_{n-1}$. Throughout, we assume that  $ [A:B]^{-1}_0 = \tau$. We refer to \cite{BGS2022} for a quick introduction to the $C^*$-index theory.
	
	Next we recall few definitions and lemmas which we need in the sequel.
	
	\begin{lmma}[\cite{BakshiVedlattice}]\label{pushdown}
		Given any $x_1\in A_1$, there exists a unique element $x_0\in A$ for which the relation $x_1e_1=x_0e_1$, holds. This element is explicitly given by $x_0=\tau^{-1} E_1(x_1e_1)$.
	\end{lmma}
	
	For each $n \geq 0$, consider the relative commutants of $B$ in $A_n$ as $B^\prime\cap A_n:= \{x\in A_n: xb=bx\mbox{ for all }b\in B\}$. These are known to be finite-dimensional, as established in \cite{Watataniindex}. On each relative commutant $B'\cap A_n$, one can define a consistent `Markov type trace', using the minimal conditional expectations (\cite{BakshiVedlattice}[Proposition $2.21$]). Specifically, for each $n\geq 0$, the map $\mathrm{tr}_n: B^{\prime}\cap A_n\to\mathbb{C}$ given by $\mathrm{tr}_n=(E_0\circ E_1\circ\cdots\circ E_n)_{|_{B^{\prime}\cap A_n}}$ defines a faithful tracial state on $B^{\prime}\cap A_n$. 
	
	\begin{ppsn}[\cite{BakshiVedlattice}]\label{m1}
		For each $n \geq 0$, there exists a faithful tracial state $\mathrm{tr}_n$ on the relative commutant $B^{\prime}\cap A_n$ satisfying the following:
		\begin{equation}
		\mathrm{tr}_n(xe_n)={\tau}\mathrm{tr}_{n-1}(x) \,\,\text{ for all } x\in B^{\prime}\cap A_{n-1}
		\end{equation}
	and for all $n\geq 1$ we have $\restr{\tr_n}{B'\cap A_{n-1}}=\tr_{n-1}$.
	\end{ppsn}
	
	For the sake of notational convenience, we shall occasionally omit the subscript and simply write $\tr$ in place of $\tr_n$. The following results will be particularly useful in what follows.
	
	\begin{lmma}[\cite{BakshiVedlattice}]\label{f2} Suppose $\{\lambda_i:1\leq i\leq n\}\subset A$ is a quasi-basis for the minimal conditional expectation $E_0$. Then, the conditional expectation from $B^{\prime}\cap A_n$ onto $A^{\prime}\cap A_n$ that preserves the trace $\mathrm{tr}$ is given by the following: $$E^{B^{\prime}\cap A_n}_{A^{\prime}\cap A_n}(x)=\tau\sum_{i}\lambda_ix\lambda^*_i\,,\quad x \in B'\cap A_n.$$
	\end{lmma}
	Let $\{\lambda_i:1\leq i\leq n\}\subset A$ be a quasi-basis associated with the minimal conditional expectation $E_0$. By proceeding similarly to the proof of [Lemma 2.23]\cite{BakshiVedlattice}, it is clear that the $\mathrm{tr}$-preserving conditional expectation from $B^{\prime}\cap A_n$ onto $A^{\prime}_1\cap A_n$ is given by $E^{B^{\prime}\cap A_n}_{A^{\prime}_1\cap A_n}(x)=\tau\sum_{i,j}\lambda_ie_1\lambda_j\,x\,\lambda^*_je_1\lambda^*_i\,,\quad x \in B'\cap A_n.$
	As noted in \cite{Watataniindex}, the following relation holds:
	\begin{equation}\label{vip}
	\sum_i \lambda_ie_1\lambda^*_i=1.
	\end{equation}
	
	\begin{crlre}[\cite{BakshiVedlattice}]\label{E-e}
		Using the notation introduced in $\Cref{f2}$, we obtain $E^{B^{\prime}\cap A_1}_{A^{\prime}\cap A_1}(e_1)=\tau.$
	\end{crlre}
	\begin{lmma}[\cite{KajiwaraWatatani}] \label{kajiwara watatani}
     Let $B\subset A$ be a unital inclusion of simple $C^*$-algebras and let $E_0$ denote the minimal conditional expectation from $B$ onto $A$, then for all $x \in B' \cap A $ and all $a \in A$, we have $E_0(xa)=E_0(ax)$ .
 \end{lmma}

  \subsection{Entropy}
 Extending the definition of Connes and St{\o}rmer \cite{CS}, Pimsner and Popa \cite{PP} defined the relative entropy $H(M|N)$ for arbitary  von Neumann subalgebras $M,N \subset Q$ as follows:
 
 \begin{dfn}[\cite{CS}]
 Let ($Q,\tau$) be a finite von Neumann algebra, and let $M,N \subset Q$ be von Neumann subalgebras. Let $\gamma = \{ w_j \in Q_+ : \sum\limits_{j=0}^{n}w_j =1\}$ be a finite partion of unity then
    $$H_\gamma(M|N) := \sum\limits_{j=0}^{n}(\tau \circ \eta E_{{N}}(w_j)- \tau \circ \eta E_{{M}}(w_j)) .$$
The Connes and St{\o}rmer relative entropy $H(M, N)$ is defined as $H({M}|{N}):= \sup_\gamma H_\gamma({M}|{N})$.
\end{dfn}
 
 Connes and St{\o}rmer developed the entropy for automorphisms of finite von Neumann algebras in \cite{CS}. Choda extended this definition of entropy to $*$-endomorphism of finite von Neumann algebras in \cite{MCJOT}. It is possible to define the entropy for $*$-endomorphism of general $C^*$-algebra \cite{CNT,C,CH}. For more details we refer the reader to \cite{NS}. Now, we briefly discuss the entropy of $*$-endomorphism of finite von Neumann algebras. Let $Q$ be a von Neumann algebra and $\tr$ be a normal trace on $Q$. Let $Q_1, \cdots, Q_n$ be finite dimensional subalgebras of $Q$, and let $E_R$ denote the $\tr$-preserving conditional expectation from $Q$ onto a von Neumann subalgebra $R$. Then,
\[
H_{\tr}(Q_1, \cdots, Q_n)=\sup\Big\{\sum_{i_1, \cdots, i_n} \eta(\tr(w_{i_1 \cdots i_n}))-\sum^{n}_{k=1}\sum_{i_k}\tau\big(\eta\big(E_{Q_k}(w^{(k)}_{i_k})\big)\big)\Big\},
\]
where the supremum is taken over all finite partitions of unity $1= \sum_{i_1 \in I_1, \cdots, i_n \in I_n}w_{i_1 \cdots i_n}$, with $w_{i_1 \cdots i_n}\geq 0$ in $Q$. Here, $w^{(k)}_{i_k}= \sum_{i_1, \cdots,i_{k-1}i_{k+1},\cdots,i_n} w_{i_1 \cdots i_n}$.
 \begin{dfn}[\cite{MCJOT}]
 Let $Q$ be a hyperfinite von Neumann algebra, $\tr$ is a normal tracial state on $Q$ and $\sigma$ is a $\tr$-preserving $*$-endomorphism of $P$. Then for each finite dimensional subalgebra $R$ of $Q$,
 $$H_{\tr, \sigma}(R)= \lim_{n\rightarrow \infty}\frac{1}{n}H_{\tr}(R, \sigma(R), \cdots, \sigma^{n-1}(R)). $$
 The entropy $H_{\tr}(\sigma)$ of $\sigma$ is the supremum of $H_{\tr, \sigma}(R)$, where $R$ runs through all finite dimensional subalgebras of $Q$.
 \end{dfn}
 \begin{dfn}[\cite{Stormer}]\label{generatingdfn}
 Let $Q$ be a finite von Neumann algebra equipped with a fixed faithful normal trace $\tr$, such that $\tr(1)=1$. Suppose that $Q_1 \subset Q_2 \subset Q_3 \subset \cdots \subset Q_n \subset \cdots$ is an increasing sequence of finite dimensional von Neumann algebras for which $Q$ is the S.O.T. closure of $\bigcup_j Q_j$. Then  the increasing sequence $
 (Q_n)_{n\in \bbn}$ of finite dimensional von Neumann subalgebras, is called a generating sequence for a $\tr$-invariant $*$-endomorphism $\sigma$ if the following conditions hold\,:
 \begin{enumerate}
 \item[(i)] $\sigma(Q_n) \subset Q_{n+1}, n\in \bbn$
 \item[(ii)] $H_{\tr}(\sigma)= \lim_{n\rightarrow \infty} \frac{1}{n}H_{\tr}(Q_n).$
\end{enumerate}
 $(Q_n)_{n\in \bbn}$ satisfies the commuting square condition if (i) holds, and $E_{Q_{j+1}}E_{\sigma( Q_{j+1})}= E_{\sigma(Q_{j})}$, for all $j\in \bbn$.
 \end{dfn}

\begin{thm}[\cite{Stormer}]\label{Stormerthm}
Let $Q$ be a finite von Neumann algebra with a fixed faithful normal trace $\tr$ with $\tr(1)=1$. Also, assume that $Q_1 \subset Q_2 \subset Q_3 \subset \cdots \subset Q_n \subset \cdots$ be an increasing sequence of finite dimensional von Neumann algebras for which $Q$ is the sot-closure of $\,\bigcup_j Q_j$. Suppose $\sigma$ is a $\tr$-invariant $*$-endomorphism with $H_{\tr}(\sigma)<\infty$, and $(Q_n)_{n\in \bbn}$ is a generating sequence for $\sigma$ satisfying the commuting square condition. Then, we have the following\,:
  \begin{enumerate}
  \item[(i)] $\lim_{n\rightarrow \infty}\frac{1}{n}H_{\tr}(Z(Q_n))$ exists.
  \item[(ii)] $H_{\tr}(\sigma)=\frac{1}{2}H_{\tr}(Q| \sigma(Q))+\frac{1}{2} \lim_{n\rightarrow \infty}\frac{1}{n}H_{\tr}(Z(Q_n)).$
  \end{enumerate}
Furthermore, if $Q$ is of type $I$, then $H_{\tr}(\sigma)=H_{\tr}(Q| \sigma(Q))$.
 \end{thm}
 
\section{Higher rotation operators and reflections}\label{Higher rotation operators}
	
Rotation maps on the $n$-box spaces of a subfactor with finite Jones index play a significant role in the development of the subfactor theory. These maps were highly instrumental in the formalism of Jones' \textit{planar algebras} ({\cite{Jo2}). They also appeared in Ocneanu's paragroup and  Popa's $\lambda$-lattice. Analogous to the subfactor theory, in \cite{BakshiVedlattice} the authors have provided a Fourier theory by defining Fourier transforms and rotation operators on the $n$-box spaces $B^{\prime}\cap A_n$ and $A^{\prime}\cap A_n$ for an inclusion of simple unital $C^*$-algebras $B\subset A$ with a conditional expectation of index-finite type (in the sense of Watatani). It was also shown that the rotation operators on the $2$-box spaces $B^\prime \cap A_1$ and $A^\prime \cap A_2$ are unital involutive $`*$' preserving anti-homomorphisms.  Recently, motivated by  \cite{Liunoncommutative}, Bakshi etal. (in \cite{BGS2022}) has established noncommutative Fourier theoretic inequalities on the $2$-box spaces of an inclusion of simple unital $C^*$-algebras using techniques from \cite{BakshiVedlattice}. To achieve this we have exploited the  fact that rotation maps on the $2$-box spaces are reflections. In this section we further investigate these operators. To extend the Fourier theoretic inequalities for the $n$-box spaces for $n> 2$, the major obstacle is the fact that the rotation operators $\rho^{\pm}_{n-1}$ are no longer reflections; that is they are not involutions for $n>2$. To see this, just observe that $(\rho^+_2)^2(e_2e_1)=e_1 e_2$. Similarly,  $\rho^{-}_2$ is also not an involution as can be readily seen from \Cref{rho+relnrho-}.} We obtain a correct notion of reflection on the $(2n+2)$-box spaces and show in the next section that using reflection we can obtain various  Fourier theoretic inequalities on the $n$-box spaces. Furthermore, we discuss canonical shifts on the $n$-box spaces which play a significant role in establishing  the relationship between Connes-St{\o}rmer entropy and Watatani index as expounded in the forthcoming section.

\smallskip

\noindent \textbf{Notation:} Throughout this section, we consider an inclusion of simple unital $C^*$-algebras $B\subset A$ with a conditional expectation $E$ of index-finite type. Let $\{\lambda_i:1\leq i\leq n\}\subset A$ be a quasi-basis for the minimal conditional expectation $E_0$. We also put $v^{(k)}_n= e_n\cdots e_k$ for $1\leq k \leq n$ and  we simply write $v_n$ instead of $v^{(1)}_n.$

 For a subfactor $N\subset M$ with $[M:N]<\infty$, the theory of Fourier transform and rotations on the relative commutants (first introduced by Ocneanu \cite{O}) had been given in great detail in \cite{Bi94} (see also \cite{Bi4,BiJo}). For recent works in this direction, see \cite{LiuFourier}. Analogously, first recall the Fourier transform and rotation operators on the relative commutants of an inclusion of simple $C^*$-algebras as defined in \cite{BakshiVedlattice}.
 
	\begin{dfn}[Fourier transform and its inverse]\label{fourier}
		For each $n\geq 0$, the Fourier transform $\mathcal{F}_n:B^{\prime}\cap A_n\longrightarrow A^{\prime}\cap A_{n+1}$ is defined by the following,
		\begin{center}
			$\mathcal{F}_n(x)=\tau^{-\frac{(n+2)}{2}}\,E^{B^{\prime}\cap A_{n+1}}_{A^{\prime}\cap A_{n+1}}(xv_{n+1})\,.$
		\end{center}
		The inverse Fourier transform ${\mathcal{F}}^{-1}_n: A^{\prime}\cap A_{n+1}\longrightarrow B^{\prime}\cap A_n$ is defined by the following,
		\begin{center}
			$\mathcal{F}^{-1}_n(w)=\tau^{-\frac{(n+2)}{2}}\,E_{n+1}(wv^*_{n+1})\,.$
		\end{center}
	\end{dfn}
 \begin{rmrk}
     The term `inverse' in the above definition is justified by the fact that $\mathcal{F}_n \circ \mathcal{F}^{-1}_n =id_{A^\prime \cap A_{n+1}} $ and $\mathcal{F}^{-1}_n \circ \mathcal{F}_n =id_{B^\prime \cap A_{n}} $ for all $n\geq 0$ (\cite{BakshiVedlattice}[Theorem 3.5]).
 \end{rmrk}
\begin{dfn}[Rotation operators]\label{rotationoperator}
	For each $n \geq 0$, the rotation map ${\rho^+_n}^{(B\subset
		A)}:B^{\prime}\cap A_n\rightarrow B^{\prime}\cap A_n$ is defined as
	\begin{equation}\nonumber
	{\rho^+_n}^{(B\subset
		A)}(x)={\big({\mathcal{F}}^{-1}_n\big({\mathcal{F}_n(x)
		}^*\big)\big)}^*,\ x \in B'\cap A_n.
	\end{equation}
Analogous to ${\rho^+_n}^{(B\subset
		A)}$, we can define a  rotation operator ${\rho^-_n}^{(B\subset
		A)}: A^\prime \cap A_{n+1} \longrightarrow A^\prime \cap A_{n+1}$ by the following,
	\begin{center}
		${\rho^-_n}^{(B\subset
			A)}(w)=\big(\mathcal{F}_n\big({\mathcal{F}_n^{-1}(w)}^*\big)\big)^* , \ w \in A^{\prime}\cap A_{n+1}.$
	\end{center}
\end{dfn}

\noindent \textbf{Notation:}
To avoid the notational difficulty, we write ${\rho^{\pm}_n}^{(B\subset
	A)}$   simply as ${\rho^{\pm}_n}$, if the inclusion $B\subset A$ is clear from the context. Also, for $0\leq k \leq n,$ we will denote $E_{k}\circ E_{k+1}\circ \cdots \circ E_n$ by $E^{n}_{k}$.
	
\medskip

Recall, by \cite{BakshiVedlattice} the rotation operators ${\rho^{\pm}_1}$ are involutions. As mentioned earlier, this fact is no longer true in higher dimensions. In this section, we find the right notion of reflection for the $(2n+2)$-box spaces.
\smallskip

Below we record two useful results whose easy proofs are omitted. 
\begin{lmma}\label{imp2}
 For each $n\geq 1$, if we take $\{\lambda_{i}: 1\leq i \leq (n+1)\}\subset A$. Then we have 
    $$\lambda_{1} {{v}^{*}_n} \lambda_{2} {{v}^{*}_{n-1}}\cdots \lambda_{n} {{v}^{*}_1}\lambda_{n+1} = \lambda_{1}{v_1} \lambda_{2} {v_2} \lambda_{3} \cdots v_n\lambda_{n+1} .$$
\end{lmma}
\begin{lmma}\label{imp} For each $1\leq k < n,$ we have,
    $$v^{(k+1)}_n(v^{(k)}_n)^*=\tau ^{n-k}e_n.$$
    \end{lmma}

  \noindent  \textbf{Notation:} In the sequel for notational convenience we denote the number  $\tau^{\big( \binom{k}{2}-kn\big)}$ by $c^{k}_n.$ Here, we use the convention that $\binom{1}{2}=0$. It is easy to verify that $\tau^{k-n}c^{k}_n=c^{k+1}_n$.

  Now we provide an important formula for $({\rho^+_n})^k$, which will be exploited extensively in the sequel.    
   \begin{ppsn}
For $1\leq k \leq n$ and $x\in B'\cap A_n  $, we have
\begin{equation*}\label{gamma+}
    ({\rho^+_n})^k(x)=c^{k}_n\sum_{i_1, \cdots, i_k}  E^{n}_{{n-k+1}} \big(v_{n-k+1}\lambda_{i_{1}}v_{n-k+2}\lambda_{i_{2}}\cdots v_{n}\lambda_{i_{k}}x\big) v_{n-k+1}\lambda^{*}_{i_{k}}v_{n-k+2}\lambda^{*}_{i_{k-1}}\cdots v_{n}\lambda^{*}_{i_{1}}.
\end{equation*}\end{ppsn}
\begin{prf}
    It is an inductive process. According to \cite{BakshiVedlattice}[Remark 3.8], we have, for all $x\in B'\cap A_n$,
    $${\rho^+_n}(x)=\tau^{-n}\sum_{i_1}E_n\big(v_n\lambda_{i_1}x\big)v_n\lambda^{*}_{i_{1}}.$$ 
The statement holds for $k=1$. Moreover, we have
    \begin{eqnarray}\nonumber
       ({\rho^+_n})^2(x)&= & \tau^{-n}\sum_{i_1}E_n\big(v_n\lambda_{i_1}{\rho^+_n}(x)\big)v_n\lambda^{*}_{i_{1}}\\\nonumber
       &=& \tau^{-2n}\sum_{i_1, i_2} E_n\big(v_n\lambda_{i_1}E_n\big(v_n\lambda_{i_2}x\big)v_n\lambda^{*}_{i_2}\big)v_n\lambda^{*}_{i_1}.
        \end{eqnarray}
One readily obtains the following equalities.  \begin{eqnarray}\nonumber
        &&v_n\lambda_{i_1}E_n(v_n\lambda_{i_2}x)v_n\\\nonumber
        &=&e_n e_{n-1}\cdots e_2 e_1\lambda_{i_1}E_n(v_n\lambda_{i_2}x)e_ne_{n-1} \cdots e_2 e_1\\\nonumber
        &=&e_n v_{n-1}\lambda_{i_1}E_n(v_n\lambda_{i_2}x)e_nv_{n-1}\\\nonumber
        &=&E_{n-1}\big(v_{n-1}\lambda_{i_1}E_n(v_n\lambda_{i_2}x)\big)e_n v_{n-1}\\\nonumber
        &=&E_{n-1}\big(v_{n-1}\lambda_{i_1}E_n(v_n\lambda_{i_2}x)\big)v_n\\\nonumber
         &=&E_{n-1}\circ E_n(v_{n-1}\lambda_{i_1} v_n\lambda_{i_2}x)v_n.
        \end{eqnarray}
Therefore, the following equations hold true.
       \begin{eqnarray}\nonumber
           ({\rho^+_n})^2(x) &=&   \tau^{-2n}\sum_{i_1, i_2} E_n\big(E_{n-1}\circ E_n(v_{n-1}\lambda_{i_1} v_n\lambda_{i_2}x)v_n\lambda^{*}_{i_2}\big)v_n\lambda^{*}_{i_1}\\\nonumber&=&   \tau^{-2n}\sum_{i_1, i_2} E_n\big(E_{n-1}\circ E_n(v_{n-1}\lambda_{i_1} v_n\lambda_{i_2}x)e_nv_{n-1}\lambda^{*}_{i_2}\big)v_n\lambda^{*}_{i_1}\\\nonumber
&=&  \tau^{-2n+1}\sum_{i_1, i_2}E_{n-1}\circ E_n(v_{n-1}\lambda_{i_1} v_n\lambda_{i_2}x)v_{n-1}\lambda^{*}_{i_2}v_n\lambda^{*}_{i_1}\\\nonumber
         &=&  c^{2}_n\sum_{i_1, i_2} E^{n}_{{n-1}}(v_{n-1}\lambda_{i_1}v_n\lambda_{i_2}x)v_{n-1}\lambda^{*}_{i_2}v_n\lambda^{*}_{i_1}.
    \end{eqnarray}
Thus, the claim also holds for $k=2$. Now, let $n\geq 3$ and assume that the statement is true for all $k<n$. Then,
    \begin{eqnarray}\nonumber
         &&({\rho^+_{n}})^{k+1}(x)\\\nonumber&=&\tau^{-n}\sum_{i_1}E_n\big(v_n\lambda_{i_1}({\rho^+_n})^{k}(x)\big)v_n\lambda^{*}_{i_{1}}\\\nonumber
         &=& \tau^{-n}c^{k}_n\sum_{i_1,\cdots, i_{k+1}} E_n\big(v_n\lambda_{i_1} E^{n}_{{n-k+1}} (v_{n-k+1}\lambda_{i_{2}} \cdots v_{n}\lambda_{i_{k+1}}x)v_{n-k+1}\lambda^{*}_{i_{k+1}}\cdots v_{n}\lambda^{*}_{i_{2}}\big)v_n\lambda^{*}_{i_{1}}.
         \end{eqnarray}
First observe that, $v_n=e_n\cdots e_{n-k+2}e_{n-k+1}e_{n-k}\cdots e_1=v^{(n-k+2)}_n e_{n-k+1} v_{n-k}$. We also recall the identity, $e_{n-k+1}ye_{n-k+1}= E_{n-k}(y)e_{n-k+1}$, which holds whenever $y \in A_{n-k}$. Now, consider $y= v_{n-k}\lambda_{i_1} E^{n}_{{n-k+1}} (v_{n-k+1}\lambda_{i_{2}} \cdots v_{n}\lambda_{i_{k+1}}x)$. Then, by applying the above identity, we obtain:
         \begin{eqnarray}\nonumber
            && v_n\lambda_{i_1} E^{n}_{{n-k+1}} (v_{n-k+1}\lambda_{i_{2}} \cdots v_{n}\lambda_{i_{k+1}}x)v_{n-k+1}\lambda^{*}_{i_{k+1}}\\\nonumber
            &=& v^{(n-k+2)}_ne_{n-k+1}v_{n-k}\lambda_{i_1} E^{n}_{{n-k+1}} (v_{n-k+1}\lambda_{i_{2}} \cdots v_{n}\lambda_{i_{k+1}}x)e_{n-k+1}v_{n-k}\lambda^{*}_{i_{k+1}}\\\nonumber
            &=&v^{(n-k+2)}_n E_{n-k}\big(v_{n-k}\lambda_{i_1} E^{n}_{{n-k+1}} (v_{n-k+1}\lambda_{i_{2}} \cdots v_{n}\lambda_{i_{k+1}}x)\big)e_{n-k+1}v_{n-k}\lambda^{*}_{i_{k+1}}\\\nonumber
             &=&v^{(n-k+2)}_n E_{n-k}\big(v_{n-k}\lambda_{i_1} E^{n}_{{n-k+1}} (v_{n-k+1}\lambda_{i_{2}} \cdots v_{n}\lambda_{i_{k+1}}x)\big)v_{n-k+1}\lambda^{*}_{i_{k+1}}.
         \end{eqnarray}
Now, it is cleat that, $v_{n-k}\lambda_{i_1} \in A_{n-k}$. Therefore, by the bi-module property of the conditional expectation, we have:
         \begin{eqnarray}\nonumber
            && v_n\lambda_{i_1} E^{n}_{{n-k+1}} (v_{n-k+1}\lambda_{i_{2}} \cdots v_{n}\lambda_{i_{k+1}}x)v_{n-k+1}\lambda^{*}_{i_{k+1}}\\\nonumber
         &=&v^{(n-k+2)}_n E_{n-k}\circ E^{n}_{{n-k+1}} (v_{n-k}\lambda_{i_1}v_{n-k+1}\lambda_{i_{2}} \cdots v_{n}\lambda_{i_{k+1}}x)v_{n-k+1}\lambda^{*}_{i_{k+1}}\\\label{eq000}
              &=&v^{(n-k+2)}_n E^{n}_{{n-k}} (v_{n-k}\lambda_{i_1} \cdots v_{n}\lambda_{i_{k+1}}x)v_{n-k+1}\lambda^{*}_{i_{k+1}}.
         \end{eqnarray}
Since $e_{m+1} \in A^\prime_{m-1}$, for all $m\geq 0$, it follows that $[v^{(n-k+2)}_n , E^{n}_{{n-k}} (v_{n-k}\lambda_{i_1} \cdots v_{n}\lambda_{i_{k+1}}x)]=0$. Using this and \Cref{eq000}, we can observe that:
        \begin{eqnarray}\nonumber
        &&({\rho^+_{n}})^{k+1}(x)\\\nonumber
            &=&\tau^{-n}c^{k}_n\sum_{i_1,\cdots, i_{k+1}} E_n\big(v^{(n-k+2)}_{n}E^{n}_{{n-k}}(v_{n-k}\lambda_{i_1}\cdots v_{n}\lambda_{i_{k+1}}x)v_{n-k+1}\lambda^{*}_{i_{k+1}}\cdots v_{n}\lambda^{*}_{i_{2}}\big)v_n\lambda^{*}_{i_{1}}\\\label{form}
         &=&\tau^{-n}c^{k}_n\sum_{i_1,\cdots, i_{k+1}} E_n\big(E^{n}_{{n-k}}(v_{n-k}\lambda_{i_1}\cdots v_{n}\lambda_{i_{k+1}}x)v^{(n-k+2)}_{n}v_{n-k+1}\lambda^{*}_{i_{k+1}}\cdots v_{n}\lambda^{*}_{i_{2}}\big)v_n\lambda^{*}_{i_{1}}.
         \end{eqnarray}
         Now we have, \begin{equation*}
             v_{n-k+1}\lambda^{*}_{i_{k+1}}\cdots v_{n}\lambda^{*}_{i_{2}}=(e_{n-k+1} e_{n-k}\cdots e_1\lambda^{*}_{i_{k+1}})\cdots(e_{n-1}e_{n-2}\cdots e_1\lambda^{*}_{i_{3}}) (e_{n}e_{n-1}\cdots e_1\lambda^{*}_{i_{2}}) .
         \end{equation*}
Since $e_{m+1} \in A^\prime_{m-1}$ for all $m\geq 0$, it follows that the first term of each of parentheses commutes with all the terms of the preceding  parentheses except the first term. Consequently, we obtain:
         
         \begin{eqnarray}\nonumber
        && v_{n-k+1}\lambda^{*}_{i_{k+1}}\cdots v_{n}\lambda^{*}_{i_{2}}\\\nonumber
        &=&e_{n-k+1}\cdots e_{n-1}e_{n}( e_{n-k}\cdots e_1\lambda^{*}_{i_{k+1}})\cdots(e_{n-2}\cdots e_1\lambda^{*}_{i_{3}}) (e_{n-1}\cdots e_1\lambda^{*}_{i_{2}})\\
       \label{999}
         &=&(v^{(n-k+1)}_n)^*v_{n-k}\lambda^{*}_{i_{k+1}}\cdots v_{n-1} \lambda^{*}_{i_{2}}.
         \end{eqnarray}
Consequently, utilizing \Cref{imp}, we deduce:

      \begin{eqnarray}\nonumber
          v^{(n-k+2)}_{n}v_{n-k+1}\lambda^{*}_{i_{k+1}}\cdots v_{n}\lambda^{*}_{i_{2}}&=&  v^{(n-k+2)}_{n} (v^{(n-k+1)}_n)^*v_{n-k}\lambda^{*}_{i_{k+1}}\cdots v_{n-1}\lambda^{*}_{i_{2}}\\\nonumber
&=& \tau^{k-1} e_n v_{n-k}\lambda^{*}_{i_{k+1}}\cdots v_{n-1}\lambda^{*}_{i_{2}}.
      \end{eqnarray}
In conclusion, by using \Cref{form}, we arrive at:
         \begin{eqnarray}\nonumber
         &&({\rho^+_{n}})^{k+1}(x)\\\nonumber
        &=&\tau^{k-(n+1)}c^{k}_n\sum_{i_1,\cdots, i_{k+1}} E_n\big(E^{n}_{{n-k}}(v_{n-k}\lambda_{i_1}\cdots v_{n}\lambda_{i_{k+1}}x)e_n v_{n-k}\lambda^{*}_{i_{k+1}}\cdots v_{n-1}\lambda^{*}_{i_{2}}\big)v_n\lambda^{*}_{i_{1}} \\\nonumber 
         &=&c^{k+1}_{n}\sum_{i_1,\cdots, i_{k+1}} E^{n}_{{n-k}}(v_{n-k}\lambda_{i_1}\cdots v_{n}\lambda_{i_{k+1}}x)v_{n-k}\lambda^{*}_{i_{k+1}}\cdots v_{n-1}\lambda^{*}_{i_{2}}v_n\lambda^{*}_{i_{1}}\\\nonumber
      &=&c^{k+1}_{n}\sum_{i_1,\cdots, i_{k+1}} E^{n}_{{n-k}}(v_{n-k}\lambda_{i_1}\cdots v_{n}\lambda_{i_{k+1}}x )v_{n-k}\lambda^{*}_{i_{k+1}}\cdots v_n\lambda^{*}_{i_{1}}.
 \end{eqnarray}
Consequently, the induction argument establishes the claim.
 \qed
\end{prf}

The following theorem, which is crucial for the subsequent analysis in this paper within the setting of inclusion of simple $C^*$-algebras, was established in the context of type $II_1$ factors in \cite{BiJo2}.

\begin{thm}\label{high rho+}
For $n\geq1$, we have $({\rho^+_n})^{n+1}=id_{B^\prime \cap A_n}$. 
\end{thm}
\begin{prf}
     The case  $n=1$, has been previously established in \cite{BakshiVedlattice}[Lemma 3.12]. For $n\geq 2$, by invoking \Cref{gamma+}, we obtain the following for every $x\in B'\cap A_n  $:
 \begin{equation*}
    ({\rho^+_n})^n(x)=c^{n}_n\sum_{i_1, \cdots, i_n}  E^{n}_{{1}} (v_{1}\lambda_{i_{1}} \cdots v_{n}\lambda_{i_{n}}x)v_{1}\lambda^{*}_{i_{n}} \cdots v_{n}\lambda^{*}_{i_{1}}.
\end{equation*}
As a result, we arrive at:

\begin{eqnarray}\nonumber
&&({\rho^+_{n}})^{n+1}(x)\\\nonumber&=&\tau^{-n}\sum_{i_1}E_n\big(v_n\lambda_{i_1}({\rho^+_n})^{n}(x)\big)v_n\lambda^{*}_{i_{1}}\\\nonumber
         &=& \tau^{-n}c^{n}_n\sum_{i_1,\cdots, i_{n+1}} E_n\big(v_n\lambda_{i_1} E^{n}_{{1}} (v_{1}\lambda_{i_{2}} \cdots v_{n}\lambda_{i_{n+1}}x)v_{1}\lambda^{*}_{i_{n+1}}\cdots v_{n}\lambda^{*}_{i_{2}}\big)v_n\lambda^{*}_{i_{1}}.
         \end{eqnarray}
It is straightforward to see that
         \begin{eqnarray}\nonumber
           v_n\lambda_{i_1} E^{n}_{{1}} (v_{1}\lambda_{i_{2}} \cdots v_{n}\lambda_{i_{n+1}}x)v_{1}\lambda^{*}_{i_{n+1}}
           &=& e_n\cdots e_2e_1\lambda_{i_1} E^{n}_{{1}} (v_{1}\lambda_{i_{2}} \cdots v_{n}\lambda_{i_{n+1}}x)e_{1}\lambda^{*}_{i_{n+1}} \\\nonumber
            &=& e_n\cdots e_2E_0\big(\lambda_{i_1} E^{n}_{{1}} (v_{1}\lambda_{i_{2}} \cdots v_{n}\lambda_{i_{n+1}}x)\big)e_{1}\lambda^{*}_{i_{n+1}}\\\nonumber
             &=& e_n\cdots e_2E_0\big( E^{n}_{{1}} (\lambda_{i_1}v_{1}\lambda_{i_{2}} \cdots v_{n}\lambda_{i_{n+1}}x)\big)e_{1}\lambda^{*}_{i_{n+1}}\\\nonumber
             &=& v^{(2)}_n E^{n}_{{0}} (\lambda_{i_1}v_{1}\lambda_{i_{2}} \cdots v_{n}\lambda_{i_{n+1}}x)v_{1}\lambda^{*}_{i_{n+1}}.
             \end{eqnarray}
Since $v^{(2)}_n \in A^\prime \subset B^\prime$, it follows that:
         
         \begin{eqnarray}\nonumber
         &&({\rho^+_{n}})^{n+1}(x)\\\nonumber
          &=& \tau^{-n}c^{n}_n\sum_{i_1,\cdots, i_{n+1}} E_n\big(v^{(2)}_nE^{n}_{0}(\lambda_{i_1}v_{1}\lambda_{i_{2}} \cdots v_{n}\lambda_{i_{n+1}}x) v_{1}\lambda^{*}_{i_{n+1}}\cdots v_{n}\lambda^{*}_{i_{2}}\big)v_n\lambda^{*}_{i_{1}}\\\label{3.4}
     &=& \tau^{-n}c^{n}_n\sum_{i_1,\cdots, i_{n+1}} E_n\big(E^{n}_{0}(\lambda_{i_1}v_{1}\lambda_{i_{2}} \cdots v_{n}\lambda_{i_{n+1}}x)v^{(2)}_n v_{1}\lambda^{*}_{i_{n+1}}\cdots v_{n}\lambda^{*}_{i_{2}}\big)v_n\lambda^{*}_{i_{1}}.
      \end{eqnarray}
By applying the same line of argument, as it was in \Cref{gamma+}, we observe that:
      \begin{eqnarray}\nonumber
           v_{1}\lambda^{*}_{i_{n+1}}\cdots v_{n}\lambda^{*}_{i_{2}}
           &=& (e_{1}\lambda^{*}_{i_{n+1}})(e_2e_{1}\lambda^{*}_{i_{n}})\cdots (e_{n-1}e_{n-2}\cdots e_1\lambda^{*}_{i_{3}}) (e_{n}e_{n-1}\cdots e_1\lambda^{*}_{i_{2}})\\\nonumber
                      &=& e_{1}e_2\cdots e_{n-1}e_{n} (\lambda^{*}_{i_{n+1}})(e_{1}\lambda^{*}_{i_{n}})\cdots (e_{n-2}\cdots e_1\lambda^{*}_{i_{3}}) (e_{n-1}\cdots e_1\lambda^{*}_{i_{2}})\\\nonumber
                       &=& v^*_{n} \lambda^{*}_{i_{n+1}}v_{1}\lambda^{*}_{i_{n}}\cdots 
 v_{n-2}\lambda^{*}_{i_{3}} v_{n-1} \lambda^{*}_{i_{2}}
 \end{eqnarray}  
Consequently, applying\Cref{3.4}, we get
      \begin{eqnarray}\nonumber
      &&({\rho^+_{n}})^{n+1}(x)\\\nonumber
         &=& \tau^{-n}c^{n}_n\sum_{i_1,\cdots, i_{n+1}} E_n\big(E^{n}_{0}(\lambda_{i_1}v_{1}\lambda_{i_{2}} \cdots v_{n}\lambda_{i_{n+1}}x)v^{(2)}_n v^*_n \lambda^{*}_{i_{n+1}} v_{1}\lambda^{*}_{i_{n}}\cdots v_{n-1}\lambda^{*}_{i_{2}}\big)v_n\lambda^{*}_{i_{1}}\\\nonumber
          &=& \tau^{-1}c^{n}_n\sum_{i_1,\cdots, i_{n+1}} E_n\big(E^{n}_{0}(\lambda_{i_1}v_{1}\lambda_{i_{2}} \cdots v_{n}\lambda_{i_{n+1}}x)e_n\lambda^{*}_{i_{n+1}} v_{1}\lambda^{*}_{i_{n}}\cdots v_{n-1}\lambda^{*}_{i_{2}}\big)v_n\lambda^{*}_{i_{1}}\quad(\text{by \Cref{imp}} )\\\nonumber
           &=& c^{n}_n \sum_{i_1,\cdots, i_{n+1}} E^{n}_{0}(\lambda_{i_1}v_{1}\lambda_{i_{2}} \cdots v_{n}\lambda_{i_{n+1}}x)\lambda^{*}_{i_{n+1}}v_{1}\lambda^{*}_{i_{n}}\cdots v_{n-1}\lambda^{*}_{i_{2}}v_n\lambda^{*}_{i_{1}} \\\nonumber
           &=&c^{n}_n \sum_{i_1,\cdots, i_{n+1}} E^{n}_{0}(\lambda_{i_{1}} {{v}^{*}_n} \lambda_{i_{2}}{{v}^{*}_{n-1}} \cdots \lambda_{i_{n}} {{v}^{*}_1} \lambda_{i_{n+1}}x)(\lambda_{i_{1}} {{v}^{*}_n} \lambda_{i_{2}}{{v}^{*}_{n-1}} \cdots \lambda_{i_{n}} {{v}^{*}_1} \lambda_{i_{n+1}})^{*}\quad (\text{by \Cref{imp2}  } ) \\\nonumber
            &=& c^{n}_n \sum_{i_1,\cdots, i_{n+1}} E^{n}_{0}(x\,\lambda_{i_{1}} {{v}^{*}_n} \lambda_{i_{2}}{{v}^{*}_{n-1}} \cdots \lambda_{i_{n}} {{v}^{*}_1} \lambda_{i_{n+1}})(\lambda_{i_{1}} {{v}^{*}_n} \lambda_{i_{2}}{{v}^{*}_{n-1}} \cdots \lambda_{i_{n}} {{v}^{*}_1} \lambda_{i_{n+1}})^{*}.\quad(\text{by \Cref{kajiwara watatani}})
\end{eqnarray}
By \cite{BakshiVedlattice}[Proposition 2.17]\footnote{It is a correct place to mention that, there was a typo in writing the quasi-basis of $E_0\circ E_1\circ \cdots E_{n-1}\circ E_n$ in \cite{BakshiVedlattice}, there the correct quasi-basis will be $ \{ \tau^{-\frac{n(n+1)}{4}}\lambda_{i_{n+1}} (e_1e_2\cdots e_{n-1}e_n)\lambda_{i_{n}} (e_1e_2\cdots e_{n-2}e_{n-1}) \cdots \lambda_{i_{3}}(e_1e_2) \lambda_{i_{2}}e_1 \lambda_{i_1}: 1\leq \lambda_{i_1}, \cdots, \lambda_{i_{n+1}}\leq m\} $.}, we know that $ \{ \tau^{-\frac{n(n+1)}{4}}\lambda_{i_{1}} {{v}^{*}_n} \cdots \lambda_{i_{n}} {{v}^{*}_1} \lambda_{i_{n+1}}: (\lambda_{i_1}, \cdots, \lambda_{i_{n+1}}) \in I^{n+1} \} $ is a quasi-basis for the conditional expectation $E^n_0$ from $A_n$ onto $B$. Thus we finally get $({\rho^+_{n}})^{n+1}(x)=x$.
\qed
\end{prf}

Next we show similar type of results for ${\rho^-_n}$. To begin with we prove a useful result.
\begin{ppsn}\label{rho+relnrho-}
  For each $n\geq 1$,
    ${\rho^{-}_n}=i\circ {\rho^{+}_n}^{(A\subset A_1)}\circ i$,
    where $i$ is the conjugate map from a $C^*$-algebra to itself.
\end{ppsn}
\begin{prf}
   Using  \Cref{rotationoperator}, it is straightforward to verify that: \begin{eqnarray}\label{eqnrho1}
       {\rho^{-}_{n}}(w)=\tau^{-(n+1)}\sum_{i_{1}}\lambda_{i_{1}}v^{*}_{n+1} E_{n+1}(w \,v^{*}_{n+1}\lambda^{*}_{i_{1}}).
  \end{eqnarray} 
and the following holds for all $x\in B'\cap A_n$, 
    $${\rho^+_n}(x)=\tau^{-n}\sum_{i_1}E_n\big(v_n\lambda_{i_1}x\big)v_n\lambda^{*}_{i_{1}}.$$
Let $\tilde{\lambda}_i:= \tau^{-\frac{1}{2}}\lambda_ie_1$. Then \{$\tilde{\lambda}_i: i \in I \}$ constitutes a quasi-basis for $E_1$. Therefore we have for $w \in A' \cap A_{2n+2}$, 
    \begin{eqnarray}\nonumber
        &&{\rho^+_n}^{(A\subset A_1)}(w^*)=\tau^{-n}\sum_{i_1}E_{n+1}\big(v^{(2)}_{n+1}\tilde{\lambda}_{i_1}\,w^{*}\big)v^{(2)}_{n+1}\tilde{\lambda}^{*}_{i_{1}}.
    \end{eqnarray}
Note that for $m\geq 2$, $v^{(2)}_{m}\tilde{\lambda}_{i_{}}=  \tau^{-\frac{1}{2}}v^{(2)}_{m}\lambda_i e_1= \tau^{-\frac{1}{2}}\lambda_i v^{(2)}_{m} e_1=\tau^{-\frac{1}{2}}\lambda_i v_{m}$ and 
 $v^{(2)}_{m}{\tilde{\lambda}}^{*}_{i}  = \tau^{-\frac{1}{2}}v^{(2)}_{m}e_1\lambda^{*}_{i}=\tau^{-\frac{1}{2}}v_{m}\lambda^{*}_{i}.$ Thus we have,
 \begin{eqnarray}\label{eqnrho2}
        &&{\rho^+_n}^{(A\subset A_1)}(w^*)=\tau^{-(n+1)}\sum_{i_1}E_{n+1}\big(\lambda_{i_1}v_{n+1}\,w^{*}\big)v_{n+1}\lambda^{*}_{i_{1}}.
    \end{eqnarray}
From \Cref{eqnrho1} and \Cref{eqnrho2}, it is straightforward to observe that the following identity holds: ${\rho^{-}_{n}}(w)=\big({\rho^+_n}^{(A\subset A_1)}(w^*)\big)^*.$ This completes the proof. \qed
\end{prf}

	 \begin{ppsn}\label{high rho -}
  For each $1\leq k \leq n$ and $w\in A'\cap A_{n+1}$, we have
     \begin{eqnarray}\nonumber
      ({\rho^-_n})^k(w)=c^{k}_{n+1}\sum_{i_1,\cdots, i_k}\lambda_{i_1}v^{*}_{n+1} \lambda_{i_2}v^{*}_{n}\cdots  \lambda_{i_k}v^{*}_{n-k+2}\, E^{{n+1}}_{{n-k+2}}(w\,v^{*}_{n+1}\lambda^{*}_{i_k}v^{*}_{n}\lambda^{*}_{i_{k-1}} \cdots v^{*}_{n-k+2}\lambda^{*}_{i_1} ).
     \end{eqnarray}
      \end{ppsn}
 \begin{prf}
Using \Cref{rho+relnrho-}, we get $({\rho^-_n})^k$ as the $n$-times composition of $i\circ {\rho^{+}_n}^{(A\subset A_1)}\circ i$. Also since $i^2= id$, we have
\begin{eqnarray}\nonumber
    &&({\rho^-_n})^k(w)\\\nonumber
    &=&i\circ {\rho^{+}_n}^{(A\subset A_1)}\circ i \circ i\circ {\rho^{+}_n}^{(A\subset A_1)}\circ i\circ \cdots \circ i\circ {\rho^{+}_n}^{(A\subset A_1)}\circ i(w)\\\nonumber
    &=&i\circ \big({\rho^{+}_n}^{(A\subset A_1)}\big)^k\circ i(w)\\\nonumber
   &=& \bigg(c^{k}_n\sum_{i_1, \cdots, i_k}  E^{n+1}_{{n-k+2}} \big(v^{(2)}_{n-k+2}\tilde{\lambda}_{i_{1}}v^{(2)}_{n-k+3}\tilde{\lambda}_{i_{2}}\cdots v^{(2)}_{n+1}\tilde{\lambda}_{i_{k}}\, w^*\big) v^{(2)}_{n-k+2}\tilde{\lambda}^{*}_{i_{k}}v^{(2)}_{n-k+3}\tilde{\lambda}^{*}_{i_{k-1}}\cdots v^{(2)}_{n+1}\tilde{\lambda}^{*}_{i_{1}}\bigg)^*.
\end{eqnarray}
Here \{$\tilde{\lambda}_i= \tau^{-\frac{1}{2}}\lambda_ie_1: i \in I \}$ is the quasi-basis of $E_1$. Now for $m\geq 2$, $v^{(2)}_{m}\tilde{\lambda}_{i_{}}=  \tau^{-\frac{1}{2}}v^{(2)}_{m}\lambda_i e_1= \tau^{-\frac{1}{2}}\lambda_i v^{(2)}_{m} e_1=\tau^{-\frac{1}{2}}\lambda_i v_{m}$ and 
 $v^{(2)}_{m}{\tilde{\lambda}}^{*}_{i}  = \tau^{-\frac{1}{2}}v^{(2)}_{m}e_1\lambda^{*}_{i}=\tau^{-\frac{1}{2}}v_{m}\lambda^{*}_{i}.$ Thus we have,
 \begin{eqnarray}\nonumber
    &&({\rho^-_n})^k(w)\\\nonumber
   &=&\bigg( c^{k}_n\tau^{-k}\sum_{i_1, \cdots, i_k}  E^{n+1}_{{n-k+2}} \big(\lambda_{i_{1}}v_{n-k+2}\lambda_{i_{2}}v_{n-k+3}\cdots \lambda_{i_{k}}v_{n+1}\, w^*\big) v_{n-k+2}\lambda^{*}_{i_{k}}v_{n-k+3}\lambda^{*}_{i_{k-1}}\cdots v_{n+1}\lambda^{*}_{i_{1}}\bigg)^*\\\nonumber
   &=& c^{k}_{n+1}\sum_{i_1,\cdots, i_k}\lambda_{i_1}v^{*}_{n+1} \cdots \lambda_{i_{k-1}}v^{*}_{n-k+3} \lambda_{i_k}v^{*}_{n-k+2}\, E^{{n+1}}_{{n-k+2}}(w\,v^{*}_{n+1}\lambda^{*}_{i_k}\cdots v^{*}_{n-k+3}\lambda^{*}_{i_{2}} v^{*}_{n-k+2}\lambda^{*}_{i_1} ).
\end{eqnarray}
\qed
 \end{prf}

 \begin{thm}\label{rotationmaporder}
	For $n\geq1$, we have $({\rho^-_n})^{n+1}=id_{A^\prime \cap A_{n+1}}$.
\end{thm}
\begin{prf}
It is a consequence of \Cref{rho+relnrho-} and \Cref{high rho+}. Composing $i\circ {\rho^{+}_n}^{(A\subset A_1)}\circ i$ $(n+1)$-times and using $i^2= id$, we have for $w\in A'\cap A_{n+1}$
     \begin{eqnarray}\nonumber
      &&({\rho^-_n})^{n+1}(w)\\\nonumber
          &=&i\circ \big({\rho^{+}_n}^{(A\subset A_1)}\big)^{(n+1)}\circ i(w)\\\nonumber
     &=&w.
      \end{eqnarray}
 The last equality follows since $\big({\rho^{+}_n}^{A\subset A_1}\big)^{(n+1)}=id_{A^\prime \cap A_{n+1}}$ by \Cref{high rho+}.
      \qed
             \end{prf}

\begin{rmrk}\label{see}

The  rotations are neither $`*$'-preserving nor anti-homomorphisms. One can easily  check that $\rho^+_2(e_2e_1)=\tau$, which is not equal to $\rho^+_2(e_1)\rho^+_2(e_2)=e_1e_2$. Furthermore $\big(\rho^+_2(e_1e_2)\big)^*=e_1e_2$, which is not equal to $\rho^+_2(e_2e_1)=\tau$. Thanks to \Cref{rho+relnrho-}, we have $\rho^{-}_2$ is also neither $`*$' preserving nor anti-homomorphisms.
\end{rmrk}

 \subsection{Reflection operators}
 In view of \Cref{high rho+} and \Cref{rotationmaporder} we introduce a new class of operators which are the right candidates for the `reflections' on the $(2n+2)$-box spaces.
	 \begin{dfn}[Reflection operators]\label{reflection operator}
	 For $n\geq 0$, define the reflection operators ${(r^+_{2n+1})}^{B\subset A}: B^\prime \cap A_{2n+1} \to B^\prime \cap A_{2n+1}$ by $$(r^+_{2n+1})^{B\subset A}=(\rho^+_{2n+1})^{n+1}.$$
 Analogous to $(r^+_{2n+1})^{B\subset A}$, we can define a reflection operator $(r^-_{2n+1})^{B\subset A}: A^\prime \cap A_{2n+2} \to A^\prime \cap A_{2n+2}$ by $$(r^-_{2n+1})^{B\subset A}=(\rho^-_{2n+1})^{n+1}.$$
	\end{dfn}

For notational convenience, we will omit the superscript $B\subset A$, and simply write $r^{\pm}_{2n+1}$ in place of $(r^{\pm}_{2n+1})^{B\subset A}$, whenever the inclusion is understood. Following the notation given in  \cite{BGS2022}, for $n=0$ we have $r^+_1=\rho^+_1=\rho_+$. We already know from \cite{BakshiVedlattice} that $r^+_1$ is a unital involutive $*$ preserving anti-homomorphisms. This is also true for $n>0$ also and we will prove it shortly. First we prove few useful results. 	 
  
  \begin{lmma}\label{imp3}
      If  $e_{[-1,2n+1]}$ denotes the Jones projection corresponding to the iterated basic construction $B \subset A_n \subset A_{2n+1}$, then we have $$e_{[-1,2n+1]} \lambda_{i_{1}}v^{*}_{n} \cdots \lambda_{i_{n}}v^{*}_1\lambda_{i_{n+1}} =\tau^{-\frac{n(n+1)}{2}}v_{n+1}\lambda_{i_{1}} v_{n+2}\lambda_{i_{2}} \cdots v_{2n+1}\lambda_{i_{n+1}}.$$
  \end{lmma}
  \begin{prf}
  Using \cite{PP2}[Theorem 2.6], it follows that for all $n\geq 1$:
  $$e_{[-1,2n+1]}= \tau^{-\frac{n(n+1)}{2}} (e_{n+1}e_n \cdots e_1)(e_{n+2}e_{n+1}\cdots e_2) \cdots(e_{2n+1}e_{2n} \cdots e_{n+1}).$$
Therefore, invoking \Cref{imp2}, we get
      \begin{eqnarray}\nonumber
          &&e_{[-1,2n+1]}  \lambda_{i_{1}}v^{*}_{n}  \cdots \lambda_{i_{n}}v^{*}_1 \lambda_{i_{n+1}} \\\label{iterated}
        &=&\tau^{-\frac{n(n+1)}{2}} (e_{n+1}e_n \cdots e_1)(e_{n+2}e_{n+1}\cdots e_2) \cdots(e_{2n+1}e_{2n} \cdots e_{n+1})\lambda_{i_1} {v_1} \lambda_{i_{2}} \cdots  v_n\lambda_{i_{n+1}}.
        \end{eqnarray}
Since for all $m\geq 0$, $e_{m+1} \in A^\prime_{m-1}$, we can observe that $\lambda_{i_1}$ commutes with all the preceding terms except the terms in the first parentheses. Thus, we get
          \begin{eqnarray}\nonumber
              &&(e_{n+1}e_n \cdots e_1)(e_{n+2}e_{n+1}\cdots e_2) \cdots(e_{2n+1}e_{2n} \cdots e_{n+1})\lambda_{i_1} {v_1} \lambda_{i_{2}} \cdots  v_n\lambda_{i_{n+1}}\\\nonumber
             &=& (e_{n+1}e_n \cdots e_1)\lambda_{i_1}(e_{n+2}e_{n+1}\cdots e_2) \cdots(e_{2n+1}e_{2n} \cdots e_{n+1}) {v_1} \lambda_{i_{2}} \cdots  v_n\lambda_{i_{n+1}}.
          \end{eqnarray}
 Now notice that $v_1 \lambda_{i_2}$ commutes with all the terms except the terms in the second parentheses and so, we get the following equality.
          \begin{eqnarray}\nonumber
           && (e_{n+1}e_n \cdots e_1)\lambda_{i_1}(e_{n+2}e_{n+1}\cdots e_2) \cdots(e_{2n+1}e_{2n} \cdots e_{n+1}) {v_1} \lambda_{i_{2}} \cdots  v_n\lambda_{i_{n+1}}   \\\nonumber
           &=& (e_{n+1}e_n \cdots e_1)\lambda_{i_1}(e_{n+2}e_{n+1}\cdots e_2){v_1} \lambda_{i_{2}} (e_{n+3}e_{n+1}\cdots e_3) \cdots(e_{2n+1}e_{2n} \cdots e_{n+1}) v_2\lambda_{i_3}\cdots  v_n\lambda_{i_{n+1}}.
          \end{eqnarray}
          Continuing in this way we obtain,
          \begin{eqnarray}\nonumber
              &&(e_{n+1}e_n \cdots e_1)(e_{n+2}e_{n+1}\cdots e_2) \cdots(e_{2n+1}e_{2n} \cdots e_{n+1})\lambda_{i_1} {v_1} \lambda_{i_{2}} \cdots  v_n\lambda_{i_{n+1}}\\\nonumber
               &=& (e_{n+1}e_n \cdots e_1) \lambda_{i_1}(e_{n+2}e_{n+1}\cdots e_2){v_1} \lambda_{i_{2}}\cdots(e_{2n+1}e_{2n} \cdots e_{n+1})v_n\lambda_{i_{n+1}}\\ \nonumber
          &=& v_{n+1}\lambda_{i_{1}}v_{n+2}\lambda_{i_{2}}\cdots v_{2n+1}\lambda_{i_{n+1}}.
              \end{eqnarray}
Therefore, by \Cref{iterated}, it follows that:
\begin{eqnarray}\nonumber
   && e_{[-1,2n+1]}  \lambda_{i_{1}}v^{*}_{n}  \cdots 
 \lambda_{i_{n}}v^{*}_1 \lambda_{i_{n+1}} \\\nonumber
        &=&\tau^{-\frac{n(n+1)}{2}} v_{n+1}\lambda_{i_{1}}v_{n+2}\lambda_{i_{2}}\cdots v_{2n+1}\lambda_{i_{n+1}}.
\end{eqnarray} This completes the proof. \qed
  \end{prf}
  
  \begin{ppsn}\label{highrotation}
  For each $n\geq 0$ we have,
      $r^+_{2n+1}= (r^+_{1})^{B \subset A_n}$. 
      \end{ppsn}
\begin{prf}
For $n=0$, it is obvious. Let $n \geq1$. Thanks to \Cref{gamma+}, we can observe that, 
    \begin{eqnarray}\nonumber
      && r^+_{2n+1}(x)\\\nonumber
      &=& ({\rho^+_{2n+1}})^{n+1}(x)\\ \label{eqn for r+}
      &=& c^{n+1}_{2n+1}\sum_{i_1, \cdots, i_{n+1}}  E^{{2n+1}}_{{n+1}} (v_{n+1}\lambda_{i_{1}}\cdots v_{2n+1}\lambda_{i_{n+1}}x)v_{n+1}\lambda^{*}_{i_{n+1}} \cdots v_{2n+1}\lambda^{*}_{i_{1}}.
    \end{eqnarray}
By \cite{BakshiVedlattice}, for all $x \in B'\cap A_1$ we have,$$(r^+_{1})^{B \subset A}(x) = \tau^{-1}\sum_{i} E_1\big(e_1\lambda_ix\big)e_1\lambda^{*}_i.$$
Recall, $ \{\tau^{-\frac{n(n+1)}{4}}\lambda_{i_{1}} {{v}^{*}_n} \cdots \lambda_{i_{n}} {{v}^{*}_1} \lambda_{i_{n+1}}: (\lambda_{i_1}, \cdots, \lambda_{i_{n+1}}) \in I^{n+1} \} $ is a quasi-basis for the conditional expectation $E^n_0$. It is easy to observe the following,

\begin{eqnarray}\nonumber
   && (r^+_{1})^{B \subset A_n}(x)\\ \nonumber
   &=& \tau^{-\frac{(n+2)(n+1)}{2}} \sum_{i_1, \cdots, i_{n+1}}  E^{{2n+1}}_{{n+1}}(e_{[-1,2n+1]} \lambda_{i_{1}}v^{*}_{n}\cdots \lambda_{i_{n}}v^{*}_1 \lambda_{i_{n+1}}x)e_{[-1,2n+1]}\lambda^{*}_{i_{n+1}}v_1\lambda^{*}_{i_{n}}\cdots v_n\lambda^{*}_{i_{1}}\\\nonumber
    &=& c^{n+1}_{2n+1}\sum_{i_1, \cdots, i_{n+1}}  E^{{2n+1}}_{{n+1}}(v_{n+1}\lambda_{i_{1}} \cdots v_{2n+1}\lambda_{i_{n+1}}x)v_{n+1}\lambda^{*}_{i_{n+1}}\cdots v_{2n+1}\lambda^{*}_{i_{1}}.\quad (\text{by \Cref{imp3}  } )
    \end{eqnarray}
Thus from \Cref{eqn for r+}, we get $(r^+_{1})^{B \subset A_n}(x)=r^+_{2n+1}(x). $
\qed
\end{prf}
\begin{thm}\label{antihomo}
	 	For $n\geq 0$, the map $r^+_{2n+1}: B^\prime \cap A_{2n+1} \to B^\prime \cap A_{2n+1}$ is a unital involutive $*$-preserving anti-homomorphism.
	 \end{thm} 
  \begin{prf}
Using \Cref{highrotation} we get $r^+_{2n+1}= (r^+_{1})^{B \subset A_n}.$
Thus, from \cite{BakshiVedlattice}[Lemma 3.12 $\&$ Lemma 3.15], we can easily see that $r^+_{2n+1}: B^\prime \cap A_{2n+1} \to B^\prime \cap A_{2n+1}$ is an involutive $*$ - preserving anti-homomorphism. Clearly $(r^+_{1})^{B \subset A}$ is unital and so, $r^+_{2n+1}$ too.
\qed
 \end{prf}
\begin{rmrk}
 The reflection operator $r^+_{2n+1}$ is involutive can also be seen using \Cref{high rho+}.
\end{rmrk}

	 \begin{lmma}\label{reflectionmaplemma}
	 	We have $\,r^{-}_{2n+1}=\mathcal{F}_{2n+1}\circ r^{+}_{2n+1}\circ \mathcal{F}_{2n+1}^{-1}$. In other words, the following diagram commutes.
   \begin{center}

   \begin{tikzcd}
{B'\cap A_{2n+1}} \arrow[r, "\mathcal{F}_{2n+1}"] \arrow[d, "r^{+}_{2n+1}" ' black]
& {A'\cap A_{2n+2}} \arrow[d, "r^{-}_{2n+1}" black] \\
{B'\cap A_{2n+1}}  \arrow[r, black, "\mathcal{F}_{2n+1}" black]
& |[black, rotate=0]| {A'\cap A_{2n+2}}
\end{tikzcd}
\end{center}
	 \end{lmma}
 
   \begin{prf}
      For $n=0$, the result has already been established in \cite{BGS2022}[Lemma 3.2]. Now, let $n>0$. Using \Cref{reflection operator} and \Cref{rotationoperator}, we derive the following:\begin{eqnarray}\nonumber
     r^{+}_{2n+1}(x)&=& (\rho^+_{2n+1})^{n+1}(x)
  =\rho^+_{2n+1}\big((\rho^+_{2n+1})^{n}(x)\big)
=\big(\mathcal{F}_{2n+1}^{-1}\big(\big(\mathcal{F}_{2n+1}\big((\rho^+_{2n+1})^{n}(x)\big)\big)^{*}\big)\big)^{*}.
 \end{eqnarray}
By \Cref{antihomo}, we know that the reflection operator $r^{+}_{2n+1}$ is  $*$ - preserving. Hence we obtain the following:
 \begin{eqnarray}\nonumber
  r^{+}_{2n+1}(x^{*})= \big(r^{+}_{2n+1}(x)\big)^{*}
= \mathcal{F}_{2n+1}^{-1}\big(\big(\mathcal{F}_{2n+1}\big((\rho^+_{2n+1})^{n}(x)\big)\big)^{*}\big).
\end{eqnarray}
Consequently, we deduce that, 
\begin{eqnarray}\label{eqnnn}
    \mathcal{F}_{2n+1}\big( r^{+}_{2n+1}(x^{*})\big)=\big(\mathcal{F}_{2n+1}\big((\rho^+_{2n+1})^{n}(x)\big)\big)^{*}.
\end{eqnarray}
From \Cref{gamma+}, it follows that
\begin{eqnarray}\nonumber
    (\rho^+_{2n+1})^{n}(x)&=&  c^{n}_{2n+1}\sum_{i_1, \cdots, i_{n}}  E^{{2n+1}}_{{n+2}} (v_{n+2}\lambda_{i_{1}} \cdots v_{2n+1}\lambda_{i_{n}}x)v_{n+2}\lambda^{*}_{i_{n}}\cdots v_{2n+1}\lambda^{*}_{i_{1}}.
\end{eqnarray}
Therefore, by \Cref{fourier}, we note that
\begin{eqnarray}\nonumber
   && \mathcal{F}_{2n+1}\big((\rho^+_{2n+1})^{n}(x)\big)\\\nonumber
   &=& \tau^{-\frac{(2n+3)}{2}}E^{B^{\prime}\cap A_{2n+2}}_{A^{\prime}\cap A_{2n+2}}\big((\rho^+_{2n+1})^{n}(x)v_{2n+2}\big)
   \\\nonumber
   &=& \tau^{-\frac{(2n+1)}{2}}\sum_{\lambda_{i_1}}{\lambda_{i_1}}(\rho^+_{2n+1})^{n}(x)v_{2n+2}{\lambda^{*}_{i_1}}\\\nonumber
   &=&  c^{n}_{2n+1}\tau^{-\frac{(2n+1)}{2}}\sum_{i_1, \cdots, i_{n+1}} {\lambda_{i_1}} E^{{2n+1}}_{{n+2}} (v_{n+2}\lambda_{i_{2}} \cdots v_{2n+1}\lambda_{i_{n+1}} x) v_{n+2}\lambda^{*}_{i_{n+1}} \cdots v_{2n+1}\lambda^{*}_{i_{2}}v_{2n+2}{\lambda^{*}_{i_1}}.
\end{eqnarray}
Applying \Cref{eqnnn}, we deduce that:
\begin{eqnarray}\nonumber
   && \mathcal{F}_{2n+1}\big( r^{+}_{2n+1}(x^{*})\big)\\\nonumber
   &=&c^{n}_{2n+1}\tau^{-\frac{(2n+1)}{2}}\sum_{i_1, \cdots, i_{n+1}} {\lambda_{i_1}}v^{*}_{2n+2}{\lambda_{i_2}}v^{*}_{2n+1}\cdots {\lambda_{i_{n+1}}}v^{*}_{n+2} E^{{2n+1}}_{{n+2}}(x^{*}\,\lambda^{*}_{i_{n+1}}v^{*}_{2n+1} \cdots \lambda^{*}_{i_{2}}v^{*}_{n+2})\lambda^{*}_{i_1}\\\label{eqn2}
    &=&c^{n}_{2n+1}\tau^{-\frac{(2n+1)}{2}}\sum_{i_1, \cdots, i_{n+1}} {\lambda_{i_1}}v^{*}_{2n+2} \cdots {\lambda_{i_{n+1}}}v^{*}_{n+2} E^{{2n+1}}_{{n+2}}(x^{*}\, \lambda^{*}_{i_{n+1}}v^{*}_{2n+1}\cdots \lambda^{*}_{i_{2}}v^{*}_{n+2}\lambda^{*}_{i_1}).
\end{eqnarray}
By applying \Cref{reflection operator} and \Cref{high rho -}, we obtain
\begin{eqnarray}\nonumber
&&r^{-}_{2n+1}\big(\mathcal{F}_{2n+1}(x^{*})\big)\\\label{r-(f(x*))}
&=&c^{n+1}_{2n+2}\sum_{i_1,\cdots, i_{n+1}} \lambda_{i_1}v^{*}_{2n+2} \cdots  \lambda_{i_{n+1}}v^{*}_{n+2} E^{{2n+2}}_{{n+2}} (\mathcal{F}_{2n+1}(x^{*}) v^{*}_{2n+2}\lambda^{*}_{i_{n+1}} \cdots v^{*}_{n+2}\lambda^{*}_{i_1}).
\end{eqnarray}
From \Cref{fourier}, it follows that
\begin{eqnarray}\nonumber
&&   \mathcal{F}_{2n+1}(x^{*}) v^{*}_{2n+2}\\\nonumber
&=& \tau^{-\frac{(2n+1)}{2}}\sum_{\lambda_{i_{n+2}}}{\lambda_{i_{n+2}}}x^{*} \,v_{2n+2}{\lambda^{*}_{i_{n+2}}}v^{*}_{2n+2}\\\nonumber
&=& \tau^{-\frac{(2n+1)}{2}}\tau^{2n+1}\sum_{\lambda_{i_{n+2}}}{\lambda_{i_{n+2}}}x^{*}\,E_{0}({\lambda^{*}_{i_{n+2}}})e_{2n+2} \quad (\text{by \cite{BakshiVedlattice}[Equation 3.1]})\\\nonumber
&=& \tau^{\frac{(2n+1)}{2}}\sum_{\lambda_{i_{n+2}}}{\lambda_{i_{n+2}}}E_{0}({\lambda^{*}_{i_{n+2}}})x^{*}\,e_{2n+2}\\\nonumber
&=& \tau^{\frac{(2n+1)}{2}} x^{*}\,e_{2n+2} .
\end{eqnarray}
Consequently, it follows that \begin{eqnarray}\nonumber
     && E^{{2n+2}}_{{n+2}} (\mathcal{F}_{2n+1}(x^{*}) v^{*}_{2n+2}\lambda^{*}_{i_{n+1}} \cdots v^{*}_{n+2}\lambda^{*}_{i_1})\\\nonumber
     &=& E^{{2n+2}}_{{n+2}} (\mathcal{F}_{2n+1}(x^{*}) v^{*}_{2n+2}\lambda^{*}_{i_{n+1}} v^{*}_{2n+1}\lambda^{*}_{i_{n}} \cdots v^{*}_{n+2}\lambda^{*}_{i_1})\\\nonumber
     &=& \tau^{\frac{(2n+1)}{2}} 
 E^{{2n+2}}_{{n+2}} ( x^{*}\,e_{2n+2} \lambda^{*}_{i_{n+1}} v^{*}_{2n+1}\lambda^{*}_{i_{n}} \cdots v^{*}_{n+2}\lambda^{*}_{i_1})\\\nonumber
  &=& \tau^{\frac{(2n+3)}{2}}
 E^{{2n+1}}_{{n+2}} ( x^{*}\, \lambda^{*}_{i_{n+1}} v^{*}_{2n+1}\lambda^{*}_{i_{n}} \cdots v^{*}_{n+2}\lambda^{*}_{i_1})
\end{eqnarray}
Finally, from \Cref{r-(f(x*))}, we obtain
\begin{eqnarray}\nonumber
    &&r^{-}_{2n+1}\big(\mathcal{F}_{2n+1}(x^{*})\big)\\\nonumber
    &=&\tau^{\frac{(2n+3)}{2}} c^{n+1}_{2n+2}\sum_{i_1,\cdots, i_{n+1}} \lambda_{i_1}v^{*}_{2n+2} \cdots  \lambda_{i_{n+1}}v^{*}_{n+2} E^{{2n+1}}_{{n+2}} ( x^{*}\, \lambda^{*}_{i_{n+1}} v^{*}_{2n+1}\lambda^{*}_{i_{n}} \cdots v^{*}_{n+2}\lambda^{*}_{i_1}).
\end{eqnarray}
Thus, comparing the above with \Cref{eqn2}, we have, for all $ x \in B'\cap A_{2n+1}$, that
\begin{eqnarray}\nonumber
   \mathcal{F}_{2n+1}\big( r^{+}_{2n+1}(x^{*})\big)= r^{-}_{2n+1}\big(\mathcal{F}_{2n+1}(x^{*})\big).
\end{eqnarray}
This completes the proof.
\qed
\end{prf}

We have seen from \cite{BakshiVedlattice}[Lemma 3.17] and \cite{BGS2022}[Proposition 3.4] that if the inclusion $B\subset A$ is irreducible, then the rotation operators $\rho^+_1$ and $\rho^-_1$ are trace preserving. In this paper we generalize this by proving that $r^+_{2n+1}$ and $r^-_{2n+1}$ are also trace preserving for any inclusion $B\subset A$ (not necessarily irreducible). To prove this result first we need the following lemma.	\begin{lmma}\label{tracelemma}
		For any $n\geq 1$,  we have $$\tr_n(xe_1)=
		\tau~\tr_n(x), ~\forall x \in A^{\prime}\cap A_n .$$

	\end{lmma}
	\begin{prf}
		The proof follows easily by using the trace-preserving conditional expectation $E^{B^\prime \cap A_n}_{A^ \prime \cap A_n}$, once we observe the following:
		$$\tr(e_1x)=\tr(E^{B^\prime \cap A_n}_{A^ \prime \cap A_n}(e_1x))= \tr(E^{B^\prime \cap A_n}_{A^ \prime \cap A_n}(e_1)x)=\tau\tr(x).$$
  \qed
	\end{prf}
	
	\smallskip
\begin{thm}\label{r+ trace preserving}
		For all $n\geq 0$, the reflection operator $r^+_{2n+1}$ is trace preserving.
	\end{thm}
\begin{prf}
	We will first do it for $r^+_{1}$.
	Let $x \in B^\prime \cap A_1$. Then,
	\begin{eqnarray}\nonumber
	\tr(r^+_{1}(x))&=& \tr(\tau^{-1}\sum_i E_1(e_1\lambda_ix)e_1\lambda_i^*) \\\nonumber
	&=& \tau^{-1}\sum_i E_0\circ E_1(E_1(e_1\lambda_ix)e_1\lambda_i^*)\\\nonumber
	&=& \sum_iE_0\circ E_1(e_1\lambda_ix\lambda_i^*)\\\nonumber&=& \tau^{-1}E_0\circ E_1(e_1E^{B^\prime \cap A_1}_{A^\prime \cap A_1}(x))\\\nonumber
	&=& \tr(E^{B^\prime \cap A_1}_{A^\prime \cap A_1}(x))\quad\quad (\text{by  \Cref{tracelemma}})\\\nonumber
	&=&\tr(x).	
	\end{eqnarray}
By \Cref{highrotation}, we have $r^+_{2n+1}=(r^+_{1})^{B \subset A_n}$ and since $r^+_{1}$ is trace preserving, we can easily see that all $r^+_{2n+1}$ are trace preserving.
 \qed
 \end{prf}
\begin{ppsn}\label{reln r- r A A_1} For each $n\geq 0$,
    $r^-_{2n+1}=(r^+_{2n+1})^{A\subset A_1}$.
\end{ppsn}
 \begin{prf}
 Using \Cref{rho+relnrho-} we can observe that,
 
     \begin{eqnarray}\nonumber
      r^-_{2n+1}(w)&=&({\rho^-_{2n+1}})^{n+1}(w)\\\nonumber
          &=&i\circ \big(({\rho^{+}_{2n+1}})^{A\subset A_1}\big)^{(n+1)}\circ i(w)\\\nonumber
          &=&i\circ (r^+_{2n+1})^{A\subset A_1}\circ i(w)\\\nonumber
     &=&(r^+_{2n+1})^{A\subset A_1}(w).
      \end{eqnarray}
      The last equality follows since we know that the reflection operator $(r^+_{2n+1})^{A\subset A_1}$ is $*$- preserving by \Cref{antihomo}.
    \qed
 \end{prf}

 \begin{thm}\label{r- every property}
      For $n\geq 0$, the reflection operator $r^-_{2n+1}: A^\prime \cap A_{2n+2} \to A^\prime \cap A_{2n+2}$, is a unital, involutive, $*$-preserving anti-homomorphism. It also preserves the trace.
 \end{thm}
 \begin{prf}
     Using \Cref{reln r- r A A_1}, \Cref{antihomo} and \Cref{r+ trace preserving} the proof is straightforward. \qed
 \end{prf}

  \begin{rmrk}
 The above fact that the reflection operator $r^-_{2n+1}$, is involutive can also be proved using \Cref{rotationmaporder}.
\end{rmrk}

The following result will be used later crucially to improve the constant of the  Young's inequality as given in \cite{BGS2022}.
\begin{thm}\label{pnormequality}
	For $x \in B^\prime\cap A_{2n+1}$ and $w \in A^\prime \cap A_{2n+2}$, we have the following
	\begin{center}
		$\|x\|_p=\|r^+_{2n+1}(x)\|_p\quad\quad\text{and}\quad \quad\|w\|_p =\|r^-_{2n+1}(w)\|_p\,$
	\end{center}
where $1\leq p \leq \infty$.
\end{thm}
\begin{prf}
Let us do the case $p=2$.	Let $x \in B^\prime\cap A_{2n+1}$.  Then,
	\begin{eqnarray}\nonumber
	\|r^+_{2n+1}(x)\|_2^2 &=& \tr(|r^+_{2n+1}(x)|^2)\\\nonumber
	&=& \tr(r^+_{2n+1}(x^*)r^+_{2n+1}(x))\\\nonumber
	&=& \tr(r^+_{2n+1}(xx^*))\quad \quad
    (\text{by \Cref{antihomo}} )\\\nonumber
	&=& \tr(xx^*)\quad \quad
    (\text{by \Cref{r+ trace preserving}} )\\\nonumber
	&=& \tr(|x|^2)\\\nonumber
	&=& \|x\|_2^2.
	\end{eqnarray}
We know that the reflection operator $r^-_{2n+1}$ is also an anti-homomorphism and trace preserving by \Cref{r- every property}. Now doing same as above for the reflection operator $r^-_{2n+1}$, we will get
	$\|r^-_{2n+1}(w)\|_2=\|w\|_2 $. Once the result is established for $p=2$, it follows  that it is also true for all $1\leq p \leq \infty$, thanks to the interpolation technique.
\end{prf}

\subsection{Canonical shift}\label{Canonical shift}
Motivated by the classification of subfactors of the hyperfinite $II_1$ factor $R$ with finite Jones index, Ocneanu had introduced a certain $*$-endomorphism which he called the canonical shift on the $n$-box spaces of the subfactor. This resembles Longo's canonical endomorphism for infinite factors which plays a significant role in the index theory for infinite factors. Analogous to the subfactor theory, we  introduce canonical shift for an inclusion of simple $C^*$-algebras with a conditional  expectation of index-finite type, as a vast generalization of canonical shift for both type $II$ and type $III$ subfactors and examine its key properties which will be instrumental in establishing a relationship between Watatani index and the Connes-St{\o}rmer relative entropy of the shift in the final section. A key difference is that unlike subfactors with finite Jones index, we do not have modular conjugation operators in our case. We exploit quasi-basis and minimal conditional expectation and use the reflection operators as the key technical tool to establish this theory. To begin, we discuss certain shift operators (see \cite{O} for $II_1$ subfactor case) on the $n$-box spaces.

We begin with the following lemma, which will be useful in our subsequent analysis. The straightforward proof is omitted.
\begin{lmma}\label{shift 10}
For every $n\geq 1$, we have $$v_{n+2}\lambda_{i_{1}}\cdots v_{2n+1}\lambda_{i_{n}} \,v_{2n}\lambda_{i_{n+1}}v_{2n+1}  \lambda_{i_{n+2}}=\tau^{2n}\lambda_{i_{1}}e_{1}\lambda_{i_{2}} \,v_{n+2}\lambda_{i_{3}}v_{n+3}  \lambda_{i_{4}}\cdots v_{2n+1}  \lambda_{i_{n+2}} .$$
\end{lmma}
The following two results are well known in the context of $II_1$ factors (for instance see \cite{David}) and will be instrumental in defining the canonical shift in our setting.

\begin{ppsn}\label{shift 1}
    For every $n\geq 0$ and $x \in B' \cap A_{2n+1}$, we have $$r^+_{2n+3}(x)= \tau^{-(2n+2)}\sum_{i_1} \lambda_{i_1} v^*_{2n+2}\,r^{+}_{2n+1}(x) \,v_{2n+3}\lambda^*_{i_1}.$$
 Moreover, $r^+_{2n+3}(x) \in A^{'}_1 \cap A_{2n+3}$.
\end{ppsn}	
\begin{prf}
Let us first treat the case when $n=0$. For any $x \in B' \cap A_{1}$, applying \Cref{eqn for r+} gives the following expression.
\begin{eqnarray}\nonumber
r^+_{3}(x)
      &=& c^{2}_{3}\sum_{i_1,i_{2}}  E_2 \circ E_{{3}}(v_{2}\lambda_{i_{1}}v_{3}\lambda_{i_{2}}x)v_{2}\lambda^{*}_{i_{2}}v_{3}\lambda^{*}_{i_{1}}\\\nonumber
    &=& c^{2}_{3}\tau \sum_{i_1,i_{2}}  E_2 (v_{2}\lambda_{i_{1}} v_{2}\lambda_{i_{2}}x)v_{2}\lambda^{*}_{i_{2}}v_{3}\lambda^{*}_{i_{1}}\\\nonumber 
     &=& c^{2}_{3}\tau^2 \sum_{i_1,i_{2}}  E_2 (e_{2}\lambda_{i_{1}} e_{1}\lambda_{i_{2}}x)v_{2}\lambda^{*}_{i_{2}}v_{3}\lambda^{*}_{i_{1}}\\\nonumber 
      &=& c^{2}_{3}\tau^3 \sum_{i_1,i_{2}}   \lambda_{i_{1}} e_{1}\lambda_{i_{2}}x\,v_{2}\lambda^{*}_{i_{2}}v_{3}\lambda^{*}_{i_{1}}\\\nonumber
      &=& c^{2}_{3}\tau^2 \sum_{i_1,i_{2}}   \lambda_{i_{1}} e_{1}E_1(e_1\lambda_{i_{2}}x)v_{2}\lambda^{*}_{i_{2}}v_{3}\lambda^{*}_{i_{1}}\\\nonumber
      &=& c^{2}_{3}\tau^2 \sum_{i_1,i_{2}}   \lambda_{i_{1}} e_1e_2E_1(e_1\lambda_{i_{2}}x)e_{1}\lambda^{*}_{i_{2}}v_{3}\lambda^{*}_{i_{1}}\\\nonumber
      &=& \tau^{-2} \sum_{i_1}   \lambda_{i_{1}} v^*_2\,r^+_1(x)\,v_{3}\lambda^{*}_{i_{1}}.
\end{eqnarray}
Now, let $n\geq 1$ and consider an element $x \in B' \cap A_{2n+1}$. By applying \Cref{eqn for r+} we obtain the following: 
     \begin{eqnarray}\label{imp shifttt}
      && r^+_{2n+3}(x)
      = c^{n+2}_{2n+3}\sum_{i_1, \cdots, i_{n+2}}  E^{{2n+3}}_{{n+2}} (v_{n+2}\lambda_{i_{1}}\cdots v_{2n+3}\lambda_{i_{n+2}}x)v_{n+2}\lambda^{*}_{i_{n+2}} \cdots v_{2n+3}\lambda^{*}_{i_{1}}.
    \end{eqnarray}
Consequently, we now obtain: 
    \begin{eqnarray}\nonumber
    &&E^{{2n+3}}_{{n+2}} (v_{n+2}\lambda_{i_{1}}\cdots v_{2n+3}\lambda_{i_{n+2}}x)\\\nonumber
    &=&E^{{2n+2}}_{{n+2}}\big( v_{n+2}\lambda_{i_{1}}\cdots v_{2n+2}\lambda_{i_{n+1}}\,  E_{2n+3}(v_{2n+3}\lambda_{i_{n+2}}x)\big)\\\nonumber
    &=& \tau E^{{2n+2}}_{{n+2}}( v_{n+2}\lambda_{i_{1}}\cdots v_{2n+2}\lambda_{i_{n+1}}v_{2n+2}  \lambda_{i_{n+2}}x)\\\nonumber
     &=& \tau^2 E^{{2n+2}}_{{n+2}}( v_{n+2}\lambda_{i_{1}}\cdots v_{2n+1}\lambda_{i_{n}} e_{2n+2}\,v_{2n}\lambda_{i_{n+1}}v_{2n+1}  \lambda_{i_{n+2}}x)\\\nonumber
    &=& \tau^2 E^{{2n+1}}_{{n+2}}\big( v_{n+2}\lambda_{i_{1}}\cdots v_{2n+1}\lambda_{i_{n}}\,  E_{2n+2}(e_{2n+2}\,v_{2n}\lambda_{i_{n+1}}v_{2n+1}  \lambda_{i_{n+2}}x )\big)\\\nonumber
    &=& \tau^3 E^{{2n+1}}_{{n+2}}( v_{n+2}\lambda_{i_{1}}\cdots v_{2n+1}\lambda_{i_{n}} \,v_{2n}\lambda_{i_{n+1}}v_{2n+1}  \lambda_{i_{n+2}}x)\\\nonumber
   &=& \tau^{2n+3} E^{{2n+1}}_{{n+2}}(\lambda_{i_{1}}e_{1}\lambda_{i_{2}} \,v_{n+2}\lambda_{i_{3}}v_{n+3}  \lambda_{i_{4}}\cdots v_{2n+1}  \lambda_{i_{n+2}}).
\end{eqnarray}
The last line follows from \Cref{shift 10}. Thus, from \Cref{imp shifttt}, we obtain:
 \begin{eqnarray}\nonumber
      && r^+_{2n+3}(x)\\\nonumber
      &=& c^{n+2}_{2n+3}\tau^{2n+3}\sum_{i_1, \cdots, i_{n+2}}  E^{{2n+1}}_{{n+2}} (\lambda_{i_{1}}e_{1}\lambda_{i_{2}} \,v_{n+2}\lambda_{i_{3}}v_{n+3}  \lambda_{i_{4}}\cdots v_{2n+1}  \lambda_{i_{n+2}})v_{n+2}\lambda^{*}_{i_{n+2}} \cdots v_{2n+3}\lambda^{*}_{i_{1}}.
    \end{eqnarray}
Furthermore, by applying \Cref{eqn for r+}, we have
\begin{eqnarray}\nonumber
&&\sum_{i_1} \lambda_{i_1} v^*_{2n+2}\,r^{+}_{2n+1}(x) \,v_{2n+3}\lambda^*_{i_1}\\\nonumber
&=&c^{n+1}_{2n+1}\sum_{i_1, \cdots, i_{n+2}} \lambda_{i_1} v^*_{2n+2}\, E^{{2n+1}}_{{n+1}} (v_{n+1}\lambda_{i_{2}}\cdots v_{2n+1}\lambda_{i_{n+2}}x)v_{n+1}\lambda^{*}_{i_{n+2}} \cdots v_{2n+1}\lambda^{*}_{i_{2}}v_{2n+3}\lambda^*_{i_1}\\\nonumber
&=&c^{n+1}_{2n+1}\sum_{i_1, \cdots, i_{n+2}} \lambda_{i_1} v^*_{n+1}\, E^{{2n+1}}_{{n+1}} (v_{n+1}\lambda_{i_{2}}\cdots v_{2n+1}\lambda_{i_{n+2}}x)(v^{(n+2)}_{2n+2})^*\,v_{n+1}\lambda^{*}_{i_{n+2}} \cdots v_{2n+1}\lambda^{*}_{i_{2}}v_{2n+3}\lambda^*_{i_1}\\\nonumber
&=&c^{n+1}_{2n+1}\sum_{i_1, \cdots, i_{n+2}} \lambda_{i_1} v^*_{n+1}\, E^{{2n+1}}_{{n+1}} (v_{n+1}\lambda_{i_{2}}\cdots v_{2n+1}\lambda_{i_{n+2}}x)v_{n+2}\lambda^{*}_{i_{n+2}} \cdots v_{2n+2}\lambda^{*}_{i_{2}}v_{2n+3}\lambda^*_{i_1}\\\nonumber
&=&c^{n+1}_{2n+1}\sum_{i_1, \cdots, i_{n+2}} \lambda_{i_1} v^*_{n+1}\, E_{n+1}\big(v_{n+1}\lambda_{i_{2}} \, E^{{2n+1}}_{{n+2}} (v_{n+2}\lambda_{i_{3}}\cdots v_{2n+1}\lambda_{i_{n+2}}x)\big)v_{n+2}\lambda^{*}_{i_{n+2}} \cdots v_{2n+3}\lambda^*_{i_1}\\\nonumber
&=&c^{n+1}_{2n+1}\tau \sum_{i_1, \cdots, i_{n+2}} \lambda_{i_1} v^*_{n+1}v_{n}\lambda_{i_{2}} \, E^{{2n+1}}_{{n+2}} (v_{n+2}\lambda_{i_{3}}\cdots v_{2n+1}\lambda_{i_{n+2}}x)v_{n+2}\lambda^{*}_{i_{n+2}} \cdots v_{2n+3}\lambda^*_{i_1}\\\nonumber
&=&c^{n+1}_{2n+1}\tau^{n+1} \sum_{i_1, \cdots, i_{n+2}} \lambda_{i_1} e_1\lambda_{i_{2}} \, E^{{2n+1}}_{{n+2}} (v_{n+2}\lambda_{i_{3}}\cdots v_{2n+1}\lambda_{i_{n+2}}x)v_{n+2}\lambda^{*}_{i_{n+2}} \cdots v_{2n+3}\lambda^*_{i_1}.
\end{eqnarray}
The last two equalities follow from \Cref{pushdown} and \cite{BakshiVedlattice}[Equation 3.2] respectively. Since $c^{n+1}_{2n+1}\tau^{n+1}\tau^{-(2n+2)}=c^{n+2}_{2n+3}\tau^{2n+3}$, so the first part is done. Once the first part is done, by proceeding in the same manner as in the proof of \cite{Bi94}[Theorem 2.13(ii)] (or by [lemma 2.23]\cite{BakshiVedlattice}), we get the other part.
 \qed
\end{prf}

Using \Cref{shift 1} we will get the following consequences whose easy proofs are same as in \cite{David}.
\begin{crlre}\label{shift 2}
\begin{itemize}
For every $n\geq 0$ and $x \in B' \cap A_{2n+1}$, we have 
   $$r^+_{2n+3}\circ r^+_{2n+1}(x)= \tau^{-(2n+2)}\sum_{i} \lambda_i v^*_{2n+2}\,x \,v_{2n+3}\lambda^*_i $$ \text{and} 
    $$ r^+_{2n+5}\circ r^+_{2n+3}(x)=r^+_{2n+3} \circ r^+_{2n+1}(x).$$
    \end{itemize}
 \end{crlre}

Now we set $S^{+}_n= r^+_{2n+3} \circ r^+_{2n+1}$, for each $n\geq 0$. Then using the above results and doing same as in the proof of \cite{Bi94}[Theorem 2.13], we will get the following observations which will be extremely crucial in the sequel.
\begin{thm} \label{shift operator}
    For each $n\geq 0$, we have
    \begin{enumerate}[label=(\roman*)]
        \item  $S^+_n:B' \cap A_{2n+1} \to A^{\prime}_1 \cap A_{2n+3}$ is a unital trace preserving $*$-isomorphism.
        \item $S^+_n|_{B' \cap A_{2n}} : B' \cap A_{2n} \to A^{\prime}_1 \cap A_{2n+2}$ is also a unital trace preserving $*$-isomorphism.
        \item $S^+_{n+1}|_{B'\cap A_{2n+1}}=S^+_n$.
       \item 
 $S^{+}_n(A^{\prime}_k \cap A_{2n+1})= A^{\prime}_{k+2} \cap A_{2n+3}$ and 
$S^{+}_n(A^{\prime}_k \cap A_{2n})= A^{\prime}_{k+2} \cap A_{2n+2}$.
    \end{enumerate}
\end{thm}

By \Cref{reln r- r A A_1}, we know that $r^-_{2n+1}=(r^+_{2n+1})^{A\subset A_1}$. Now if we set $S^{-}_n= r^-_{2n+3} \circ r^-_{2n+1}$, for each $n\geq 0$. Thus, we get the following results for  $S^-_{n}$.
\begin{thm} \label{shift 3}
    For each $n\geq 0$, we have
    \begin{enumerate}[label=(\roman*)]
        \item  $S^-_n:A' \cap A_{2n+2} \to A^{\prime}_2 \cap A_{2n+4}$ is a unital trace preserving $*$-isomorphism.
        \item $S^-_n|_{A' \cap A_{2n+1}} : A' \cap A_{2n+1} \to A^{\prime}_2 \cap A_{2n+3}$ is also a unital trace preserving $*$-isomorphism.
        \item $S^-_{n+1}|_{A'\cap A_{2n+2}}=S^-_n$. 
       \item
$S^{-}_n(A^{\prime}_{k+1} \cap A_{2n+2})= A^{\prime}_{k+3} \cap A_{2n+4}$ and
$S^{-}_n(A^{\prime}_k \cap A_{2n+1})= A^{\prime}_{k+2} \cap A_{2n+3}.$
    \end{enumerate} 
\end{thm}

As discussed in \Cref{preliminaries}, starting from an inclusion $B\subset A$ of simple unital $C^*$-algebras with a conditional expectation of index-finite type, we obtain a tower $\{A_j\}_j$ of simple unital $C^*$-algebras, by iterating the $C^*$-basic construction with respect to the minimal conditional expectations. This yields an increasing sequence of finite dimensional von Neumann algebras as follows:
\begin{center}
		$A'\cap A_1 \subset A' \cap A_2 \subset A'\cap A_3 \subset \cdots \subset A'\cap A_n \subset \cdots$
	\end{center}
	Recall further that the maps $\mathrm{tr}_n:=( E_1\circ E_2 \circ \cdots\circ E_n)_{|_{A^{\prime}\cap A_n}}$ define faithful tracial states on $A^{\prime}\cap A_n$. We set $P_0 = \bigcup_{n\geq 1} A'\cap A_n$. There is a natural trace on $P_0$ defined by $\tr(w)=\tr_n(w)$, where $n$ is such that $w \in A' \cap A_n$. Let $\mathcal{H}$ is the completion of $P_0$ with respect to the inner product $\langle w_1, w_2\rangle = \tr(w^*_1w_2)$. An element $w\in P_0$ when viewed as a vector in $\mathcal{H}$ will be written as $\widehat{w}$. For $w_1, w_2 \in P_0$, we set $\pi(w_1)\widehat{w}_2= \widehat{w_1w_2}$. It is easy to check that $\pi : P_0 \rightarrow B(\mathcal{H})$ is an injective $*$-homomorphism and we will write $w$ for $\pi(w)$. We denote by  $P$ the SOT closure of $P_0$ in $B(\mathcal{H})$. For $w\in P_0$, we have $\tr(w)=\big<\widehat{1}, w\widehat{1}\big>$, and we extend $\tr$ to $P$ by the same expression.

 In view of \Cref{shift 3}, the trace preserving $*$-endomorphism $\Gamma$ on $\bigcup_{n\geq 1} A'\cap A_n$ can be defined by $$\Gamma(w)= S^-_n(w), \quad \text{ for } w\in A' \cap A_{2n+2}.$$ We can extend $\Gamma$ (will denote it by the same $\Gamma$) to a trace preserving endomorphism of the von Neumann algebra $P$. Motivated by Ocneanu, we  call this endormphism $\Gamma$ the \textit{canonical shift of the inclusion $B\subset A$}. The following observations are well-known in the context of an inclusion of factors $N\subset M$, with a conditional expectation from $M$ onto $N$ (see \cite{CH}, for instance).
\begin{ppsn}\label{shiftnew}For all natural number $j$ and $k$,
    \begin{enumerate}[label=(\roman*)]

        \item $\Gamma^k(A'\cap A_j)= A^{\prime}_{2k} \cap A_{j+2k}$.
        \item $\Gamma(A'\cap A_{2j})= r^-_{2j+1}(A'\cap A_{2j})$.
        \end{enumerate}
\end{ppsn}
\begin{prf}
\begin{enumerate}[wide=0pt, listparindent=1.25em, parsep=0pt]
\item[$(i)$] It is an easy consequence of \Cref{shift 3} item $(iv)$.
\item[$(ii)$] Using \Cref{shift 3}, we have
$$\Gamma(A'\cap A_{2j})=S^-_{j-1}(A'\cap A_{2j})= A^{\prime}_2 \cap A_{2j+2}.$$
Moreover, by \Cref{reln r- r A A_1}, we know that $r^-_{2j+1}=(r^+_{2j+1})^{A\subset A_1}$. Let $w \in A' \cap A_{2j}$. Then, by \Cref{shift 1}, we have $r^-_{2j+1}(w) \in A^{\prime}_2 \cap A_{2j+2}$. Therefore, $ r^-_{2j+1}(A'\cap A_{2j})\subset A^{\prime}_2 \cap A_{2j+2}. $

Now let $\tilde{w}\in A^{\prime}_2 \cap A_{2j+2}.$ Since $S^-_{j-1}$ is an isomorphism, there exists $w\in A'\cap A_{2j}$ such that $$\tilde{w}=S^-_{j-1}(w)=r^-_{2j+1} \circ r^-_{2j-1}(w)= r^-_{2j+1}(w_1)$$ where $w_1= r^-_{2j-1}(w) \in A'\cap A_{2j}$. This shows that $A^{\prime}_2 \cap A_{2j+2} \subset r^-_{2j+1}(A'\cap A_{2j}) $. Combining this with the inclusion obtained earlier, we conclude that $r^-_{2j+1}(A'\cap A_{2j})=A^{\prime}_2 \cap A_{2j+2}=\Gamma(A'\cap A_{2j}).$
\end{enumerate}
\end{prf}

\begin{ppsn}\label{shift5} Let $E_Q$ denote the conditional expectation of $P$ onto a von Neumann subalgebra $Q$ of $A$. For $j\geq 2$, we have the relation $E_{A' \cap A_j}E_{\Gamma(A' \cap A_j)}= E_{\Gamma(A' \cap A_{j-2})}$, that is, the following is a commuting square:
 \[\begin{matrix}
A' \cap A_j&\subset & P\\
\cup & &\cup\\
 \Gamma(A' \cap A_{j-2}) &\subset & \Gamma(A' \cap A_j).
 \end{matrix}\]
    
\end{ppsn}
\begin{prf} First, observe that by applying item $(i)$ of \Cref{shiftnew}, we have $\Gamma(A' \cap A_{j-2})=A^{\prime}_2 \cap A_{j}$ and $\Gamma(A' \cap A_{j})=A^{\prime}_2 \cap A_{j+2}$. Thus we have to show that the following is a commuting square:
 \[\begin{matrix}
A' \cap A_j&\subset & P\\
\cup & &\cup\\
 A^{\prime}_2 \cap A_{j} &\subset & A^{\prime}_2 \cap A_{j+2}.
 \end{matrix}\]
To complete the argument, it suffices to show that $E_{A' \cap A_j}(A^{\prime}_2 \cap A_{j+2})=A^{\prime}_2 \cap A_{j}$. First, note that for $w \in A^{\prime}_2 \cap A_{j+2}$ and $w'\in {A^{\prime} \cap A_{j}} $, we have $$\tr\big(E^{j+2}_{j+1}(w)w')=\tr(ww').$$  Moreover, for any $w \in A^{\prime}_2 \cap A_{j+2}$ and $a \in A_2$, the following commutation relation holds: $$E^{j+2}_{j+1}(w)a=aE^{j+2}_{j+1}(w).$$
This implies that $E_{A' \cap A_j}(A^{\prime}_2 \cap A_{j+2})\subset A^{\prime}_2 \cap A_{j}$. The reverse inclusion is immediate, hence the desired equality follows.
\qed
 
\end{prf}

Now we will demonstrate that the canonical shift $\Gamma$ is a $2$-shift in the sense of \cite{Choda}. We begin by recalling the definitions of $n$-shift. 
 \begin{dfn}[\cite{Choda}]\label{nshift}
 Let $Q$ be an injective finite von Neumann algebra equipped with a faithful normal trace $\tr$ such that $\tr(1)=1$. Suppose that $Q_1 \subset Q_2 \subset Q_3 \subset \cdots \subset Q_n \subset \cdots$ is an increasing sequence of finite dimensional von Neumann algebras for which $Q$ is the SOT closure of $\bigcup_j Q_j$. For a natural number $n$, a $\tr$ preserving $*$-endomorphism $\sigma$ of $Q$ is called an $n$-shift on the tower $Q_1 \subset Q_2 \subset Q_2 \subset \cdots \subset Q_n \subset \cdots$ for $Q$, if the following conditions are satisfied:
 \begin{enumerate}[label=(\roman*)]
 \item For all $j,m \in \bbn$, the von Neumann algebra $\{Q_j, \sigma(Q_j), \cdots, \sigma^m(Q_j)\}''$ is contained in $Q_{j+nm}$.
 \item There exists a sequence $(k_j)_{j\in \bbn}$ of integers satisfying $\lim_{j\rightarrow \infty }\frac{nk_j -j}{j}=0$ such that the following properties hold, $w_1\sigma^m(w_2)=\sigma^m(w_2)w_1$ and $\tr(w_3\sigma^{lk_j}(w_1))=\tr(w_3)\tr(w_1)$, for all $l\in \bbn$ and $w_1, w_2 \in Q_j, m\in k_j\bbn$ and $w_3\in \{Q_j, \sigma^{k_j}(Q_j), \cdots, \sigma^{(l-1)k_j}(Q_j)\}''$.
 \item Let $E_R$ denote the conditional expectation of $Q$ onto a von Neumann subalgebra $R$ of $Q$. Then for all $j\in \bbn$, we have  $E_{Q_j}E_{\sigma( Q_j)}= E_{\sigma(Q_{j-n})}$, i.e. the following is a commuting square:
 \[\begin{matrix}
Q_j&\subset & Q\\
\cup & &\cup\\
 \sigma(Q_{j-n}) &\subset & \sigma( Q_j).
\end{matrix}\]
 \item For all $j\geq n$, there exists a $\tr$-preserving $*$-automorphism or anti-automorphism $\beta$ of $Q_{nj+n}$ such that $\sigma(Q_{nj})= \beta(Q_{nj})$.
 \end{enumerate}
 \end{dfn}
\begin{ppsn}\label{ppsn for 2 shift}
The canonical shift $\Gamma$ for the inclusion $B\subset A$ is a $2$-shift on the tower $(A'\cap A_j)_j$.
 \end{ppsn} 
 \begin{prf}
  Using item $(i)$ of \Cref{shiftnew}, we obtain that for each $j ,m \in \bbn$, the following identity holds: $$\Gamma^m(A'\cap A_j)= A^{\prime}_{2k} \cap A_{j+2k}.$$ From this, it follows immediately that $\{A' \cap A_j , \,\Gamma(A' \cap A_j ), \,\cdots , \,\Gamma^{m}(A' \cap A_j )\}^{\prime \prime} \subset A' \cap A_{j+2m}$.  Hence, the first condition of \Cref{nshift} is verified in this case for $n=2$.
  
  If we take $k_j= [\frac{j}{2}]+1$, then it is straightforward to verify that $$\lim_{j\rightarrow \infty }\frac{nk_j -j}{j}=0.$$ Moreover, for any $m\in k_j\bbn$ we have $\Gamma^m(A'\cap A_j)= A^{\prime}_{2m} \cap A_{j+2m}\subset A^{\prime}_j$. Hence, we obtain the following commuting relation for any $w_1, w_2 \in (A'\cap A_j)$: $$w_1\Gamma^m(w_2)=\Gamma^m(w_2)w_1.$$ It is also worth noting, by \Cref{shiftnew} item $(i)$, that we have $\Gamma^{lk_j}(A' \cap A_j)= A^{\prime}_{2lk_j}\cap A_{j+2lk_j}$ for all $l\in \bbn$. Assume that $w \in A' \cap A_j$ and $y \in \{A' \cap A_j , \,\Gamma(A' \cap A_j ), \,\cdots , \,\Gamma^{(l-1)k_j}(A' \cap A_j )\}$. Then we have $y\,\Gamma^{lk_j}(w)\in  A^{\prime}\cap A_{j+2lk_j}$, and hence we obtain the following:
\begin{eqnarray}\nonumber
\tr\big(y\,\Gamma^{lk_j}(w)\big)&=& E^{j+2lk_j}_0(y\,\Gamma^{lk_j}(w)\big)\\\nonumber
&=& E^{2lk_j}_0E^{j+2lk_j}_{2lk_j+1}(y\,\Gamma^{lk_j}(w)\big)\\\nonumber
&=& E^{2lk_j}_0\big(y\,E^{j+2lk_j}_{2lk_j+1}\big(\Gamma^{lk_j}(w)\big)\big).
\end{eqnarray}
Note that $E^{j+2lk_j}_{2lk_j+1}\big(\Gamma^{lk_j}(w)\big)$ belongs to $ A^{\prime}_{2lk_j}\cap A_{2lk_j}=\bbc$. Hence, we conclude: 
\begin{eqnarray}\nonumber
\tr\big(y\,\Gamma^{lk_j}(w)\big)&=& E^{2lk_j}_0\big(y\,E^{j+2lk_j}_{2lk_j+1}\big(\Gamma^{lk_j}(w)\big)\big)\\\nonumber
&=& E^{2lk_j}_0(y) E^{j+2lk_j}_{2lk_j+1}\big(\Gamma^{lk_j}(w)\big)\\\nonumber
&=& E^{2lk_j}_0(y) E^{j+2lk_j}_{0}\big(\Gamma^{lk_j}(w)\big)\\\nonumber
&=& \tr(y) \,\tr\big(\Gamma^{lk_j}(w)\big)\\\nonumber
&=& \tr(y) \,\tr(w).
\end{eqnarray}
As the above relation holds for any $w \in A' \cap A_j$ and any $y \in \{A' \cap A_j , \,\Gamma(A' \cap A_j ), \,\cdots , \,\Gamma^{(l-1)k_j}(A' \cap A_j )\}$, it follows that $$\tr\big(y\,\Gamma^{lk_j}(w)\big)=\tr(y) \,\tr(w)$$ for any $w \in A' \cap A_j$ and $y \in \{A' \cap A_j , \,\Gamma(A' \cap A_j ), \,\cdots , \,\Gamma^{(l-1)k_j}(A' \cap A_j )\}^{\prime \prime}$. Consequently, this establishes that the second condition in \Cref{nshift} is also fulfilled here for the case $n=2$.

By applying \Cref{shift5}, it readily follows that the third condition in \Cref{nshift} is satisfied in this case for $n=2$.
 
 Since $r^{-}_{2j+1}$ is a $\tr$ preserving $*$-anti-automorphism of ${A'\cap A_{2j+2}}$ by \Cref{r- every property}, the third condition for the \Cref{nshift} is also satisfied once we apply \Cref{shiftnew} item $(ii)$. Therefore, $\Gamma$ is a $2$-shift on the tower $(A'\cap A_j)_j$.
 \qed
 \end{prf}

\section{Convolution product and Fourier theoretic inequalities}\label{Fourier transforms revisited}

Given an inclusion $B\subset A $ of simple unital $C^*$-algebras with a conditional expectation of index-finite type, in  \cite{BakshiVedlattice} we defined a convolution product on $B^\prime \cap A_1$ and proved that the rotation operator $\rho_1$ is anti-multiplicative with respect to this new product (see also [Proposition 3.9]\cite{BGS2022}). In this section we discuss the convolution product on the $n$-box spaces $B^{\prime}\cap A_n$ and $A^\prime \cap A_{n+1}$. The key ingredient  is the higher order Fourier transform. We also obtain various Fourier theoretic inequalities generalizing \cite{BakshiVedlattice}. Furthermore, we improve the constant of the Young's inequality given in \cite{BakshiVedlattice}.
	\begin{dfn}\label{convolution}(Convolution)
		The convolution product of two elements $x$ and $y$ in $B^{\prime}\cap A_n$, denoted by $x*y$, is defined as
		\begin{center}
			$x* y=\mathcal{F}_n^{-1}\big(\mathcal{F}_n(y)\mathcal{F}_n(x)\big).$
		\end{center}
	Similarly, for any two elements $w,z \in A^\prime \cap A_{n+1}$ we define
	\begin{center}
		$w*z=\mathcal{F}_n\big(\mathcal{F}_n^{-1}(z){\mathcal{F}_n}^{-1}(w)\big).$
	\end{center}
	\end{dfn}
It was already shown in \cite{BakshiVedlattice}[Lemma 3.20] that the convolution $`*$' is associative. It was also shown in \cite{BGS2022} that the rotation operator $\rho_1^+=\rho_+$ and $\rho_1^-=\rho_-$(notations given in \cite{BGS2022}) are anti-multiplicative with respect to the convolution. But the higher rotation operator $\rho_n$ for $n\geq$ are no longer anti-multiplicative. 
However, we show that reflection operators are anti-multiplicative with respect to the convolution as seen below.
	\begin{ppsn}\label{convolutionantihomo}
		For $x,y\in B^\prime\cap A_{2n+1}$, we have $r^+_{2n+1}(x\,*\,y)=r^+_{2n+1}(y)\,*\, r^+_{2n+1}(x)$. Similarly, for $w,z\in A^{\prime}\cap A_{2n+2}$ one has $\,r^-_{2n+1}(w\, *\, z)=r^-_{2n+1}(z)\,* \,r^-_{2n+1}(w)$.
	\end{ppsn}
\begin{prf}
By \Cref{convolution} and \Cref{reflectionmaplemma}, we can easily observe that,
\begin{eqnarray}\nonumber
    r^+_{2n+1}(y)\,*\, r^+_{2n+1}(x)&=& \mathcal{F}_{2n+1}^{-1}\big(\mathcal{F}_{2n+1}(r^+_{2n+1}(x))\mathcal{F}_{2n+1}(r^+_{2n+1}(y))\big)\\\nonumber
    &=& \mathcal{F}_{2n+1}^{-1}\big(r^-_{2n+1}(\mathcal{F}_{2n+1}(x))r^-_{2n+1}(\mathcal{F}_{2n+1}(y))\big).
    \end{eqnarray}
    It is known from \Cref{r- every property}, that the reflection operator $r^-_{2n+1}$ is an anti-homomorphism. Therefore, we obtain:

    \begin{eqnarray}\nonumber
    r^+_{2n+1}(y)\,*\, r^+_{2n+1}(x) &=& \mathcal{F}_{2n+1}^{-1}\big(r^-_{2n+1}(\mathcal{F}_{2n+1}(y)\mathcal{F}_{2n+1}(x))\big).
    \end{eqnarray}
Now using once again \Cref{reflectionmaplemma} and \Cref{convolution}, we get
    \begin{eqnarray}\nonumber
    r^+_{2n+1}(y)\,*\, r^+_{2n+1}(x) &=& r^+_{2n+1}\big(\mathcal{F}_{2n+1}^{-1}(\mathcal{F}_{2n+1}(y)\mathcal{F}_{2n+1}(x))\big)\\\nonumber
      &=& r^+_{2n+1}(x* y).
\end{eqnarray}
Similarly, we can show that $\,r^-_{2n+1}(w\, *\, z)=r^-_{2n+1}(z)\,* \,r^-_{2n+1}(w)$.
	\end{prf}\qed

\begin{rmrk} By \cite{BGS2022}[Proposition 3.8], we have $\,(x* y)^*=(x^*)* (y^*)$, for any $\,x,y\in B^\prime \cap A_{1}$ and similarly, $\,(w* z)^*=(w^*)* (z^*)$ for any $w, z\in A^{\prime}\cap A_{2}$.  However, the above result does not hold for all $(2n+2)$-box spaces $B'\cap A_{2n+1}$ and $A'\cap A_{2n+2}$. For instance,  if we take $n=2$ and $x=y=e_1 \in B'\cap A_5$, then by \Cref{fourier}, we can observe that
$$\mathcal{F}_5(e_1) = \tau^{-\frac{7}{2}}E^{B^{\prime}\cap A_{6}}_{A^{\prime}\cap A_{6}}(e_1v_{6})= \tau^{-\frac{3}{2}}e_{6}e_5e_4e_3.$$
Thus, by applying \Cref{convolution} and \Cref{fourier}, we obtain:
$$e_1* e_1 ={\mathcal{F}_5}^{-1}\big(\mathcal{F}_5(e_1)\mathcal{F}_5(e_1)\big)=\tau^{-1}{\mathcal{F}_5}^{-1}(e_{6}e_4e_5e_3)=\tau^{-\frac{9}{2}} E_6(e_{6}e_4e_5e_3v^*_6)=\tau^{-\frac{1}{2}} e_4e_1e_3.$$
As a result, we get, $(e^*_1)* (e^*_1)= e_1 * e_1=\tau^{-\frac{1}{2}} e_4e_1e_3 \text{  and  }
(e_1 * e_1)^*=\tau^{-\frac{1}{2}} e_3e_1e_4.$
If the mentioned result is true then we will get  $\tau^{-\frac{1}{2}} e_4e_1e_3=\tau^{-\frac{1}{2}} e_3e_1e_4$. That implies $E^{B^{\prime}\cap A_{6}}_{A^{\prime}\cap A_{6}}(e_4e_1e_3)=E^{B^{\prime}\cap A_{6}}_{A^{\prime}\cap A_{6}}(e_3e_1e_4)$, which gives $e_4e_3=e_3e_4$. This  is not possible since $e_4e_3=e_3e_4$ implies $e_3 \in \{e_4\}'$ and so, it would imply $B=A$.
\end{rmrk}
\medskip

Instead of the inclusion of simple unital $C^*$-algebras $B\subset A$, considering the inclusion $A \subset A_1$ we obtain similarly Fourier transformation, inverse Fourier transform and the convolution product for the corresponding relative commutants. We fix the notations $\mathcal{F}^{A\subset A_1}_{n}$ and $*_1$ for the above Fourier transformation from $A'\cap A_{n+1}$ onto $A^{\prime}_1 \cap A_{n+2}$ and the corresponding convolution product of two elements in $A' \cap A_{n+1}$ respectively. Below we establish that the Fourier transforms $\mathcal{F}^{B\subset A}_{n}$ and $\mathcal{F}^{A\subset A_1}_{n}$ are closely related.

\begin{thm}\label{shift odd}
For every $n \geq 0$, we have the following:
     $$S_{[\frac{n}{2}]}\big(\mathcal{F}^{-1}_{n}(w)\big)=\big(\mathcal{F}^{A\subset A_1}_{n}(w^*)\big)^*, \quad \text{ for all } w\in A'\cap A_{n+1}.$$ In other words, if we denote the conjugate map from a  $C^*$-algebra to itself by $i$, then the following diagram commutes.

    \begin{center}
    \begin{tikzcd}
{B' \cap A_{n}} \arrow[r, "\mathcal{F}_{n}" ] \arrow[dr, "S_{[\frac{n}{2}]}" ' ]
& {A'\cap A_{n+1}} \arrow[d, "i \circ \mathcal{F}^{A\subset A_1}_{n}\circ i"]\\
& {A^{\prime}_1 \cap A_{n+2}}
\end{tikzcd}
\end{center}
\end{thm}
\begin{prf}
First, consider the case when $n$ is odd, i.e., $n=2k+1$ for some $k\geq 0$. Then using \Cref{shift operator} and \Cref{fourier}, we can easily notice that 
  \begin{eqnarray}\nonumber
    S_{[\frac{n}{2}]}\big(\mathcal{F}^{-1}_{n}(w)\big)
    &=& S_k\big(\mathcal{F}^{-1}_{2k+1}(w)\big)\\\nonumber
    &=& \tau^{-(2k+2)}\sum_{i} \lambda_i v^*_{2k+2}\, \mathcal{F}^{-1}_{2k+1}(w) \,v_{2k+3}\lambda^*_i \\\nonumber
     &=& \tau^{-(2k+2)} \tau^{-\frac{(2k+3)}{2}}\sum_{i} \lambda_i v^*_{2k+2}\, E_{2k+2}(w\, v^*_{2k+2}) \,v_{2k+3}\lambda^*_i .
    \end{eqnarray}
Since $e_{m+1} \in A^\prime_{m-1}$, for all $m\geq 0$, thus using \Cref{pushdown} and [Equation 3.2]\cite{BakshiVedlattice}, we get
    \begin{eqnarray}\nonumber
     S_{[\frac{n}{2}]}\big(\mathcal{F}^{-1}_{n}(w)\big) &=&  \tau^{-(2k+2)} \tau^{-\frac{(2k+3)}{2}}\sum_{i} \lambda_i v^*_{2k+3}\, E_{2k+2}(w\, v^*_{2k+2}) \,v_{2k+2}\lambda^*_i \\\nonumber
       &=&  \tau^{-(2k+1)} \tau^{-\frac{(2k+3)}{2}}\sum_{i} \lambda_i v^*_{2k+3}\, w\, v^*_{2k+2} \,v_{2k+1}\lambda^*_i \\\label{forAsubsetA1}
       &=&  \tau^{-\frac{(2k+3)}{2}}\sum_{i} \lambda_i v^*_{2k+3}\, w\, v_{1}\lambda^*_i. 
    \end{eqnarray} 
 Now, by applying \Cref{fourier}, we can easily observe that for all $w\in A' \cap A_{2k+2}$,
 $$\mathcal{F}^{A\subset A_1}_{n}(w^*)=\mathcal{F}^{A\subset A_1}_{2k+1}(w^*)=\tau^{-\frac{(2k+3)}{2}}\,E^{A^{\prime}\cap A_{2k+3}}_{A^{\prime}_1\cap A_{2k+3}}(w^* \,e_{2k+3}e_{2k+2} \cdots e_2).$$
Since $\{ \lambda_i : i\in I\}$ is a quasi-basis for the conditional expectation $E_0$. It is easy to check that, $\{\tau^{-\frac{1}{2}}\lambda_i e_1:  i \in I \}$ is a quasi-basis for the conditional expectation $E_1$, thus we have from \Cref{f2},

 \begin{eqnarray}\nonumber
  \mathcal{F}^{A\subset A_1}_{n}(w^*)&=&\tau^{-\frac{(2k+3)}{2}}\sum_i \lambda_i e_1 w^* \,e_{2k+3}e_{2k+2} \cdots e_2 e_1 \lambda^*_i \\\nonumber
  &=& \tau^{-\frac{(2k+3)}{2}}\sum_i \lambda_i v^*_1 w^* \,v_{2k+3} \lambda^*_i.
\end{eqnarray}
This yields:
\begin{equation}\label{AsubsetA1}
    \big(\mathcal{F}^{A\subset A_1}_{n}(w^*)\big)^*=\tau^{-\frac{(2k+3)}{2}}\sum_i \lambda_i v^*_{2k+3} w \,v_{1} \lambda^*_i.
\end{equation}
Comparing \Cref{AsubsetA1} and \Cref{forAsubsetA1}, we are done for $n=2k+1$. Now assume $n=2k$ for some $k\geq 0$. Since $e_{m+1} \in A^\prime_{m-1}$, for all $m\geq 0$, then from \Cref{shift operator}, we have for all $ x\in B' \cap A_{2k}$,
   \begin{eqnarray}\nonumber
       S_{[\frac{n}{2}]}(x)&=& S_k(x)\\\nonumber
       &=& \tau^{-(2k+2)}\sum_{i} \lambda_i v^*_{2k+2}\, x \,v_{2k+3}\lambda^*_i \\\nonumber
        &=& \tau^{-(2k+2)}\sum_{i} \lambda_i v^*_{2k+1}\, x \,e_{2k+2}e_{2k+3}e_{2k+2}\cdots e_1\lambda^*_i \\\nonumber
         &=& \tau^{-(2k+1)}\sum_{i} \lambda_i v^*_{2k+1}\, x \,v_{2k+2}\lambda^*_i.
   \end{eqnarray}
Using the same line of argument as above $n=2k+1$ case, here we will get
   \begin{eqnarray}\nonumber
       S_{[\frac{n}{2}]}\big(\mathcal{F}^{-1}_{n}(w)\big) 
       &=& \tau^{-(2k+1)}\sum_{i} \lambda_i v^*_{2k+1}\, \mathcal{F}^{-1}_{2k}(w)\,v_{2k+2}\lambda^*_i\\\nonumber
        &=&  \tau^{-\frac{(2k+2)}{2}}\sum_{i} \lambda_i v^*_{2k+2}\, w\, v_{1}\lambda^*_i.
       \end{eqnarray}
In a similar manner as in the $n=2k+1$ case, here also we can observe that 
       \begin{equation*}
  \big(\mathcal{F}^{A\subset A_1}_{n}(w^*)\big)^*=   \big(\mathcal{F}^{A\subset A_1}_{2k}(w^*)\big)^*=\tau^{-\frac{(2k+2)}{2}}\sum_i \lambda_i v^*_{2k+2} w \,v_{1} \lambda^*_i.
\end{equation*}
Hence we are done.
    \qed
\end{prf}

 Now we are ready to obtain the higher dimensional Fourier theoretic inequalities in the sense of \cite{Liunoncommutative}. Jiang, Liu and Wu provided a non-commutative version of the Hausdorff-Young inequality, the Young’s inequality and uncertainty principles for  $\mathscr{P}_{2,\pm}$ exploiting the planar algebraic formulation of the standard invariant  in \cite{Liunoncommutative}. Subsequently in \cite{ BGS2022}, we have established those inequalities for $B^{\prime}\cap A_1$ for an inclusion of simple unital $C^*$-algebras $B\subset A$. Motivated by the subfactor theory we denote the vector spaces $B'\cap A_{n-1}$ and $A'\cap A_{n}$ by $\mathscr{P}_{n,+}$ and $\mathscr{P}_{n,-}$ respectively. Our goal of this subsection is to establish various Fourier theoretic inequalities for the $n$-box spaces $\mathscr{P}_{n,\pm}$, $n\geq 2$ of  $B \subset A$. This naturally generalizes some of the results in section 4 of \cite{BGS2022}. \\
	
	\textbf{Notation}: 
		\[{\kappa}^{+}_n(B\subset A) = \text{min}\big\{ \tr(p): p\in \mathcal{P}(B^{\prime}\cap A_n)\big\},\] 
	\[{\kappa}^{-}_n (B\subset A)= \text{min}\big\{ \tr(q): q\in \mathcal{P}(A^{\prime}\cap A_{n+1})\big\}\] and 
	\[{\kappa}_n(B\subset A) =\sqrt{{\kappa}^{+}_n(B\subset A) {\kappa}^{-}_n(B\subset A)}.\]
	If there is no danger of confusion, we often write ${\kappa}^{\pm}_n$ and $\kappa_n$ in place of ${\kappa}^{\pm}_n(B\subset A)$ and ${\kappa}_n(B\subset A)$, for brevity.
	
	\medskip
	
	We now list a few important results whose easy proofs are similar to \cite{BGS2022} and  omitted.
	\begin{ppsn}\label{constant}
		For any element $x\in B^{\prime}\cap A_n$, we have the following: $$\mathcal{F}_n(x){\mathcal{F}_n(x)}^*=\tau^{-1}\,E^{B^{\prime}\cap A_{n+1}}_{A^{\prime}\cap A_{n+1}}(xe_{n+1}x^*).$$
	\end{ppsn}
	
 \begin{ppsn}\label{Fn trace preserving}
     Both $\fn $ and $\fn^{-1}$ are isometries with respect to the norm given by $\|x\|_2=(\tr(x^*x))^{\frac{1}{2}}.$
 \end{ppsn}
 We now proceed to prove the Hausdorff-Young inequality for the $n$-box spaces.  The proof strategy is similar to that of \cite{BGS2022}. Recall that $\delta=\sqrt{{[A:B]}_0}$. In the subsequent theorem, we will use the Kadison–Schwarz inequality, which states that for a positive linear functional $f$ on a $C^*$-algebra $Q$, and for $x,y\in Q$, we have:
 $$|f(y^*x)|^2 \leq f(x^*x)f(y^*y).$$

	\begin{thm}[Hausdorff-Young inequality]\label{HY}
		Let $B\subset A$ be an inclusion of simple unital $C^*$-algebras with a conditional expectation of index-finite type. \textcolor{black}{Then, we have the following\,:}
		\begin{enumerate}[label=(\roman*)]
		\item For any $x \in B^\prime \cap A_n$, we have
			$\Vert x \rVert_q \leq \lVert \mathcal{F}_n(x)\rVert_p\leq \Big(\dfrac{\delta}{\kappa_{n-1}}\Big)^{1-\frac{2}{p}}~\lVert x\rVert_q$,
		where $2 \leq p \leq \infty$ and $\frac{1}{p}+\frac{1}{q}=1$.
		\item For any $w \in A^\prime \cap A_{n+1}$, we have
			$\Vert w \rVert_q \leq \lVert \mathcal{F}^{-1}_n(w)\rVert_p\leq \Big(\dfrac{\delta}{\kappa_{n-1}}\Big)^{1-\frac{2}{p}}~\lVert w\rVert_q$, where $2 \leq p \leq \infty$ and $\frac{1}{p}+\frac{1}{q}=1$.
\end{enumerate}	
	 
	\end{thm}
	
	\begin{prf}
	\begin{enumerate}[wide=0pt, listparindent=1.25em, parsep=0pt]
	\item[$(i)$]
	Following [Theorem 4.1] \cite{BGS2022}, it sufficient to show 
	\begin{equation}\label{maintheorem}
	\,\lVert x\rVert _1\leq {\lVert \mathcal{F}_n (x)\rVert}_{\infty} \leq \dfrac{\delta}{\kappa_{n-1}} {\lVert x\rVert}_1.
	\end{equation}
	Using the Kadison-Schwarz inequality, we obtain the estimate ${\lVert x \rVert}_1 \leq {\lVert x \rVert}_{2}$. Additionally, it holds that ${\lVert x \rVert}_2\leq {\lVert x \rVert}_{\infty}$. Since $\mathcal{F}_n$ is isometry with respect to ${\lVert . \rVert}_2$, i.e., ${\lVert \mathcal{F}_n(x) \rVert}_{2}={\lVert x \rVert}_2$, we conclude that $\,\lVert x\rVert _1\leq {\lVert \mathcal{F}_n (x)\rVert}_{\infty}$. We will first prove the remaining part of the inequality for a partial isometry. Then, by applying a rank-one decomposition, we will extend the result to the general case. Similar to  \cite{BGS2022}[Lemma 4.4], we can show that for a partial isometry $v \in B^\prime \cap A_n$, the following properties hold, $E_n(v^*v)\leq \dfrac{1}{{\kappa}_{n-1}^{+}} {\lVert v\rVert}_1$ and $E^{B^{\prime}\cap A_n}_{A^{\prime}\cap A_n}(vv^*)\leq \dfrac{1}{{\kappa}_{n-1}^{-}}{\lVert v\rVert}_1$. Thus, we derive $v e_{n+1} v ^* \leq \dfrac{{\lVert v \rVert}_1}{{\kappa}_{n-1}^{+}} ~v v^*$. Furthermore, as a consequence of \Cref{constant} and the results established above, we obtain the following:
		$$\mathcal{F}_n (v){\mathcal{F}_n (v)}^*=\delta^2 \,E^{B^{\prime}\cap A_{n+1}}_{A^{\prime}\cap A_{n+1}}(ve_{n+1}v^*)\leq\dfrac{{\lVert v \rVert}_1}{{\kappa}_{n-1}^{+}}\delta^2 \,E^{B^{\prime}\cap A_{n+1}}_{A^{\prime}\cap A_{n+1}}(vv^*)\leq \dfrac{{\delta}^2{\lVert v \rVert}^2_1}{\kappa^2_{n-1}}.$$
		Applying ${\lVert . \rVert}_{\infty}$ to both sides leads us to the desired inequality.
		\item[$(ii)$]
		Using the Kadison–Schwarz inequality and proceeding as in $(i)$, we obtain the following $\Vert w \rVert_1 \leq \lVert \mathcal{F}^{-1}_n(w)\rVert_\infty $. Now, we will show that ${\lVert \mathcal{F}^{-1}_n (w)\rVert}_{\infty} \leq \dfrac{\delta}{\kappa_{n-1}} {\lVert w\rVert}_1.$ It is enough to show the inequality for a partial isometry. Now, by \Cref{maintheorem}, we have for a partial isometry $v\in A' \cap A_{n+1}$,
     $$\lVert \mathcal{F}^{A\subset A_1}_n(v^*)\rVert_\infty\leq \dfrac{\delta} {{{\kappa}}_{n-1}(A \subset A_1)} {\lVert v^*\rVert}_1.$$
 Here 
  ${{\kappa}}_{n-1}(A \subset A_1)=\sqrt{{{\kappa}}^{+}_{n-1}(A \subset A_1){{\kappa}}^{-}_{n-1}(A \subset A_1)}$. By definition ${{\kappa}}^{+}_{n-1}(A \subset A_1)={\kappa}^{-}_{n-1}$. Also by \Cref{shift operator}, we know that $A^{\prime}_1\cap A_{n+1}$ is isomorphic to $B' \cap A_{n-1}$. Hence ${{\kappa}}^{-}_{n-1}(A \subset A_1)={\kappa}^{+}_{n-1}$. So we get ${{\kappa}}_{n-1}(A \subset A_1)={\kappa}_{n-1}$ and 
  $$\lVert \mathcal{F}^{A\subset A_1}_n(v^*)\rVert_\infty\leq \dfrac{\delta} {{\kappa}_{n-1}} {\lVert v^*\rVert}_1.$$
Since $v$ is a partial isometry, one can observe that, $${\lVert v^*\rVert}_1= \tr(|v^*|)=\tr((vv^*)^{\frac{1}{2}})=\tr(vv^*)=\tr(v^*v)=\tr((v^*v)^{\frac{1}{2}})=\tr(|v|)={\lVert v\rVert}_1 .$$
As a consequence of \Cref{shift odd}, we have 
  $$S_{[\frac{n}{2}]}\big(\mathcal{F}^{-1}_{n}(w)\big)=\big(\mathcal{F}^{A\subset A_1}_{n}(w^*)\big)^*, \quad \text{ for all } w\in A'\cap A_{n+1}.$$
This leads to:
   $$\lVert \mathcal{F}^{A\subset A_1}_n(v^*)\rVert_\infty=\lVert \big(\mathcal{F}^{A\subset A_1}_n(v^*)\big)^*\rVert_\infty=\lVert S_{[\frac{n}{2}]}\big(\mathcal{F}^{-1}_{n}(v)\big)\rVert_\infty=\lVert \mathcal{F}^{-1}_{n}(v)\rVert_\infty.$$
Putting everything together, we get:
   $${\lVert \mathcal{F}^{-1}_n (v)\rVert}_{\infty} \leq \dfrac{\delta}{\kappa_{n-1}} {\lVert v\rVert}_1.$$
Thus, by applying rank-one decomposition of any $w\in A' \cap A_{n+1}$, we get
   \begin{equation}\label{-nbox imp}
       \,\lVert w\rVert _1\leq {\lVert \mathcal{F}^{-1}_n (w)\rVert}_{\infty} \leq \dfrac{\delta}{\kappa_{n-1}} {\lVert w\rVert}_1.
   \end{equation}
Since both $\mathcal{F}^{-1}_n$ and $\mathcal{F}_n$ are isometries with respect to the norm given by $\|x\|_2=(\tr(x^*x))^{\frac{1}{2}}$, then using the same kind of interpolation as in proof of \cite{BGS2022}[Theorem 4.1], we are done.\qed

		\end{enumerate}
	
	\end{prf}
	
\begin{rmrk}
 This is an appropriate place to mention that, by applying \Cref{HY}, the quintuple $(\mathscr{P}_{n,+}$, $\mathscr{P}_{n,-},\tr_{n-1}, \tr_n, d^+_n\mathcal{F}_{n-1})$ forms a von Neumann $\big(d^+_n\big)^2$-bi-algebra (see \cite{HLJ}), while the quintuple $(\mathscr{P}_{n,-}$, $\mathscr{P}_{n,+},\tr_n, \tr_{n-1}, d^-_n \mathcal{F}^{-1}_{n-1})$ forms a von Neumann $\big(d^-_n\big)^2$-bi-algebra, where $\frac{1}{d^+_n}= \sup\{{\lVert \mathcal{F}_{n-1}(x) \rVert}_{\infty}: x\in \mathscr{P}_{n,+}, {\lVert x \rVert}_1=1$\} and $\frac{1}{d^-_n}= \sup\{{\lVert \mathcal{F}^{-1}_{n-1}(w) \rVert}_{\infty}: w\in \mathscr{P}_{n,-}, {\lVert w \rVert}_1=1$\}.
	\end{rmrk}
	As a consequence of \Cref{HY}, we have the Donoho-Stark uncertainty principle and the Hirschman-Beckner uncertainty principle for the $n$-box spaces as follows (proof is similar to \cite{BGS2022,Liunoncommutative}). 
	\begin{thm}[Donoho-Stark uncertainty principle]\label{DS}
		Consider an inclusion of simple unital $C^*$-algebras $B\subset A$ with a conditional expectation of index-finite type. \textcolor{black}{Then, we obtain the following:}
		\begin{enumerate}[label=(\roman*)]
		\item 
		For any non zero $x \in B^\prime \cap A_n$, we have $
		\mathcal{S}(x)\mathcal{S}(\fn(x))\geq \dfrac{\kappa^2_{n-1}
		}{{[A:B]}_0}\,.$
		\item For any non zero $w \in A^\prime \cap A_{n+1}$, we have $
		\mathcal{S}(w)\mathcal{S}(\fn^{-1}(w))\geq \dfrac{\kappa^2_{n-1}
		}{{[A:B]}_0}\,.$
		\end{enumerate}
	\end{thm}
	\begin{thm}[Hirschman-Beckner uncertainty principle]\label{HBuncertainty}
		Let $B\subset A$ be an inclusion of simple unital $C^*$-algebras with a conditional expectation of index-finite type. \textcolor{black}{Then, we obtain the following:}
		\begin{enumerate}[label=(\roman*)]
		\item
		For any $x \in B^\prime \cap A_n$,
		$\frac{1}{2}\big(H(|\fn(x)|^2)+H(|x|^2)\big)\geq -\|x\|_2^2\,\Big(\log\bigg(\dfrac{\delta}{\kappa_{n-1}}\bigg) +\log \|x\|_2^2\Big)\,.$
		\item For any $w \in A^\prime \cap A_{n+1}$,
		$\frac{1}{2}\big(H(|\fn^{-1}(w)|^2)+H(|w|^2)\big)\geq -\|w\|_2^2\,\Big(\log\bigg(\dfrac{\delta}{\kappa_{n-1}}\bigg) +\log \|w\|_2^2\Big)\,.$

		\end{enumerate}
	\end{thm}

 As a further application, we  improve the constant of the Young's inequality for $\mathscr{P}_{2,+}$ as proved in \cite{BGS2022}. We also prove the Young's inequality  for  $\mathscr{P}_{2,-}$. First we need a Lemma. 
\begin{lmma}\label{new convolution}
    For $w,z \in A' \cap A_{n+1}$, we have,
    $ w*_1 z = (z^* * w^*)^*.$
\end{lmma}
\begin{prf}
Using \Cref{convolution}, we can easily observe that
   \begin{eqnarray}\label{new convolution 1}
    w*_1 z=   \big(\mathcal{F}^{A\subset A_1}_{n}\big)^{-1}\big( \mathcal{F}^{A\subset A_1}_{n}(z)\mathcal{F}^{A\subset A_1}_{n}(w)\big)
   \end{eqnarray} 
By applying \Cref{shift odd}, we get the following: 
   $$\mathcal{F}^{A\subset A_1}_{n}(w)=\big(S_{[\frac{n}{2}]}\big(\mathcal{F}^{-1}_{n}(w^*)\big)\big)^*=S_{[\frac{n}{2}]}\big(\big(\mathcal{F}^{-1}_{n}(w^*)\big)^*\big), \quad \text{ for all } w\in A'\cap A_{n+1}.$$
The last equality follows since $S_{[\frac{n}{2}]}$ is $*$-preseving. Also since $S_{[\frac{n}{2}]}$ is a homomorphism we get
   \begin{eqnarray}\nonumber
      \mathcal{F}^{A\subset A_1}_{n}(z)\mathcal{F}^{A\subset A_1}_{n}(w)&=& S_{[\frac{n}{2}]}\big(\big(\mathcal{F}^{-1}_{n}(z^*)\big)^*\big)S_{[\frac{n}{2}]}\big(\big(\mathcal{F}^{-1}_{n}(w^*)\big)^*\big)\\\nonumber
      &=&S_{[\frac{n}{2}]}\big(\big(\mathcal{F}^{-1}_{n}(z^*)\big)^*\big(\mathcal{F}^{-1}_{n}(w^*)\big)^*\big)\\\label{new convolution 2}
      &=&\big(S_{[\frac{n}{2}]}\big(\mathcal{F}^{-1}_{n}(w^*)\mathcal{F}^{-1}_{n}(z^*)\big)\big)^*.
      \end{eqnarray}
       Using \Cref{shift odd}, we have for $y \in A^{\prime}_1 \cap A_{n+2}$, 
     $$\big(\mathcal{F}^{A\subset A_1}_{n}\big)^{-1}(y)=i\circ \mathcal{F}_{n} \circ S_{[\frac{n}{2}]}^{-1} \circ i (y).$$
Therefore, by applying \Cref{new convolution 1} and \Cref{new convolution 2}, we ultimately obtain:
\begin{eqnarray}\nonumber
    w*_1 z&=&i\circ \mathcal{F}_{n} \circ S_{[\frac{n}{2}]}^{-1} \circ i \big(\mathcal{F}^{A\subset A_1}_{n}(z)\mathcal{F}^{A\subset A_1}_{n}(w)\big)\\\nonumber
    &=&i\circ \mathcal{F}_{n} \circ S_{[\frac{n}{2}]}^{-1} \circ i \big(\big(S_{[\frac{n}{2}]}\big(\mathcal{F}^{-1}_{n}(w^*)\mathcal{F}^{-1}_{n}(z^*)\big)\big)^*\big)\\\nonumber
     &=&i\circ \mathcal{F}_{n} \big(\mathcal{F}^{-1}_{n}(w^*)\mathcal{F}^{-1}_{n}(z^*)\big)\\\nonumber
     &=& (z^* * w^*)^*.
\end{eqnarray}
\qed
\end{prf}

We now prove Young's inequality for $\mathscr{P}_{2,-}$ and also improve the constant appearing in the Young's inequality  for $\mathscr{P}_{2,+}$ as in \cite{BGS2022}.
\begin{thm}[Young's inequality]\label{Young}
     Suppose $B\subset A $ is an inclusion of simple unital $C^*$-algebras with a conditional expectation of index-finite type. \textcolor{black}{Then, we have the following\,:}
     \begin{enumerate}[label=(\roman*)]
     \item For any $x,y \in B' \cap A_1$, we have
     $\,\lVert x * y\rVert _r \leq\dfrac{\delta}{\kappa^{+}_{0}}\lVert x \rVert _p \lVert  y\rVert _q$,
     where $1\leq p,q,r\leq \infty$ and $\frac{1}{p}+\frac{1}{q}= \frac{1}{r}+1$.
     \item For any $w,z \in A' \cap A_2$, we have
     $\,\lVert w * z\rVert _r \leq\dfrac{\delta}{\kappa^{-}_{0}}\lVert w \rVert _p \lVert  z\rVert _q$,
     where $1\leq p,q,r\leq \infty$ and $\frac{1}{p}+\frac{1}{q}= \frac{1}{r}+1$.
     \end{enumerate}
 \end{thm}
 \begin{prf}
 \begin{enumerate}[wide=0pt, listparindent=1.25em, parsep=0pt]
\item[$(i)$] This follows from \cite{BGS2022}[Theorem 4.12] and \Cref{pnormequality}.
\item[$(ii)$]
 From \Cref{new convolution} and \cite{BGS2022}[Proposition 3.8], it follows that for all $w,z \in A' \cap A_2$, we have:
     $$z*_1w=(w^* * z^*)^*= w * z.$$
     It is known that ${{\kappa}}^{+}_{0}(A \subset A_1)={\kappa}^{-}_{0}$. Hence, by \cite{BGS2022}[Lemma 4.21] and \Cref{pnormequality}, we obtain, for all $1\leq p \leq \infty$, the following:
$${\lVert w*z\rVert}_{p} ={\lVert z*_1w\rVert}_{p}\leq \dfrac{\delta}{\kappa^-_{0}} {\lVert z\rVert}_{p}{\lVert w\rVert}_1 \quad \text{ and }\quad{\lVert z*w\rVert}_{p} ={\lVert w*_1z\rVert}_{p}\leq \dfrac{\delta}{\kappa^-_{0}} {\lVert z\rVert}_{p}{\lVert w\rVert}_1 .$$
Now, using \cite{BGS2022}[Lemma 4.22], \Cref{convolutionantihomo}  and \Cref{pnormequality}, we will get for $1\leq p \leq q$ and $\frac{1}{p}+\frac{1}{q}=1$, 
$${\lVert w*z\rVert}_{\infty} \leq \dfrac{\delta}{\kappa^-_{0}} {\lVert w\rVert}_{p}{\lVert z\rVert}_q .$$
Finally doing the same kind of interpolation as in the Proof of \cite{BGS2022}[Theorem 4.12], we are done.
\qed
     \end{enumerate}
 \end{prf}
 
 We end this section with the following question.
 \smallskip
 
 \noindent\textbf{Question:} What is the Young's inequality for the higher $n$-box spaces $\mathscr{P}_{n,\pm}$?
 
 \section{Connes-St{\o}rmer entropy of the Canonical shift}\label{final section}

In this short section, we establish a relationship between Watatani index and the Connes-St{\o}rmer entropy of the canonical shift associated with an inclusion $B \subset A $ of simple unital $C^*$-algebras with a conditional expectation of index-finite type, assuming further that the inclusion is of finite depth \cite{JOPT}. Recall, a unital inclusion $B \subset A $ of $C^*$-algebras is said to be of finite depth if there is a $n\geq 1$ such that $(B' \cap A_{n-1}) e_{n} (B'\cap A_{n-1})= B'\cap A_{n}$. The smallest such $n$ is called the depth of the inclusion. In our case, it is worth noting that $\sup_n \dim Z(A' \cap A_{2n})< \infty$ and there exist a $k_0 \in \bbn$ such that $(A' \cap A_{n}) e_{n+1} (A'\cap A_{n})= A'\cap A_{n+1}$ for all $n\geq 2k_0$. Since $E_{n+1}$ is the trace preserving conditional expectation from $A'\cap A_{n+1}$ onto $A'\cap A_{n}$, then doing same as in the proof of \cite{GHJ}[Theorem 4.6.3] we get that if $G$ denotes the inclusion matrix for $A' \cap A_{2k_{0}} \subset A' \cap A_{2k_{0}+1} $ then $G$ is connected and the inclusion matrix for $A' \cap A_{2k_{0}+1} \subset A' \cap A_{2k_{0}+2}$ is $G^t$. Thereby the inclusion matrix $GG^t$ for the inclusion $A' \cap A_{2k_{0}} \subset A' \cap A_{2k_{0}+2}$ is primitive by \cite{GHJ}[Lemma 1.3.2]. Furthermore if $s^{2k_{0}}$ denotes the trace vector for the $\tr$ on $A' \cap A_{2k_0}$, then $GG^t s^{2k_0}=[A:B]_0 s^{2k_0}$. Moreover it is also easy to observe that the $GG^t$ is the inclusion matrix for the inclusions $A' \cap A_{2k} \subset A' \cap A_{2k+2}$ for all $k\geq k_0$. Recall that, as defined in \Cref{Canonical shift}, $P$ denotes the von Neumann algebra generated by $\bigcup_{n\geq 1} A'\cap A_n$. Our main objective in this subsection is to establish the following result:

\begin{thm}\label{entropy 2}
Let  $B \subset A $ be an inclusion simple unital $C^*$-algebras with a conditional expectation of index-finite type and assume that the inclusion is of finite depth. Then 
$$\frac{1}{2}H_{\tr}(P| \Gamma(P))= H_{\tr}(\Gamma)= \log[A:B]_0 .$$
\end{thm}

In proving \Cref{entropy 2}, we will need the following results.

\begin{thm}[\cite{Choda}]\label{nshiftthm}
 If $Q$, $\{Q_j:j\in \bbn\}$, $\sigma$and $\tr$ are as in \Cref{nshift}, and furthermore if $\sigma$ satisfies conditions (1) and (2) in \Cref{nshift}, then $H_{\tr}(\sigma)= \lim_{k\rightarrow \infty} \frac{H_{\tr}(Q_{nk})}{k}$.
 \end{thm}
As established earlier in \Cref{ppsn for 2 shift}, the map $\Gamma$ satisfies the conditions (1) and (2) of \Cref{nshift}. A straightforward consequence of \Cref{nshiftthm} is the following useful fact. \begin{ppsn}\label{2shift}
The canonical shift $\Gamma$ for the inclusion $B\subset A$ has the property that $H_{\tr}(\Gamma)={\lim}_{n\rightarrow \infty} \frac{H_{\tr}(A' \cap A_{2n})}{n}$.
\end{ppsn}

At this point, we recall the following result (considered as folklore), which will play a crucial role in proving \Cref{entropy 2}.

\begin{ppsn}[\cite{NS}, Proposition 10.4.9]\label{imp lma for entropy}
Let $Q=
\underset{\xrightarrow{\hspace{1em}}}{\lim}\, Q_k$ be an AF algebra for which the inclusion matrices for the inclusions $Q_k \subset Q_{k+1}$ for every $k$ are given by a fixed primitive matrix $\Lambda$. Then there exists a unique tracial state on $Q$ and
$$\lim_{n\rightarrow \infty}\frac{1}{n}H_{\tr}(Q_n)=\log\beta,$$
where $\beta$ is the Perron-Frobenious eigen value of $\Lambda$.
\end{ppsn}

\noindent\textbf{Proof of \Cref{entropy 2}: } As already discussed at the beginning of this section, for a finite-depth inclusion $B\subset A$ of simple unital $C^*$-algebras with a conditional expectation of index-finite type, there exists $k_0 \in \bbn$ such that the inclusion matrices for the inclusions $A' \cap A_{2k} \subset A' \cap A_{2k+2}$ for all $k\geq k_0$ are same, more precisely the inclusion matrix is $GG^t$, where $G$ is the inclusion matrix for $A' \cap A_{2k_{0}} \subset A' \cap A_{2k_{0}+1} $. Furthermore if $s^{2k_{0}}$ denotes the trace vector for the $\tr$ on $A' \cap A_{2k_0}$, then $GG^t s^{2k_0}=[A:B]_0 s^{2k_0}$. Observe that $s^{2k_{0}}$ is a Perron-Frobenius eigen vector of $GG^t$ with Perron-Frobenius eigen value $[A:B]_0$. Then, by applying \Cref{imp lma for entropy} to $Q_n=A'\cap A_{2n}$, and subsequently using \Cref{2shift}, we obtain $$ H_{\tr}(\Gamma)=\lim_{n\rightarrow \infty}\frac{1}{n}H_{\tr}(A' \cap A_{2n})=\log[A:B]_0.$$
Using \Cref{shiftnew} item $(i)$ and \Cref{2shift}, we observe that $\{A' \cap A_{2n}\}_n$ is a generating sequence (see  \Cref{generatingdfn}) for the canonical shift $\Gamma$. Furthermore, by applying \Cref{shift5}, we see that it satisfies the commuting square condition. 
Thanks to \Cref{Stormerthm}, $\lim_{n\rightarrow \infty}\frac{1}{n}H_{\tr}(Z(A' \cap A_{2n}))$ exists in this setting, and we obtain the following equation.
\begin{equation}\nonumber
H_{\tr}(\Gamma)=\frac{1}{2}H_{\tr}(P| \Gamma(P))+ \lim_{n\rightarrow \infty}\frac{1}{2n}H_{\tr}(Z(A' \cap A_{2n})).
\end{equation}
For finite depth inclusion $B \subset A $ we obviously have $\lim_{n\rightarrow \infty}\frac{1}{2n}H_{\tr}(Z(A' \cap A_{2n}))=0.$ Accordingly, the preceding analysis lead us to conclude the following:
 
 $$\frac{1}{2}H_{\tr}(P| \Gamma(P))= H_{\tr}(\Gamma)= \log[A:B]_0 .$$
This completes the proof.
\qed
\section*{Acknowledgements}
K.C.B acknowledges the support of INSPIRE Faculty grant DST/INSPIRE/04/2019/002754. SM wishes to thank the Indian Statistical Institute, Bangalore, where part of the work was carried out, and acknowledges support of NBHM Postdoctoral Fellowship 0204/16(20)/2022/R\&D-15506. SM also acknowledges the project grant received under PM USHA Scheme G.O number G.O.(Rt)No.239/2025/HEDN dated 22.02.2025.

	\bigskip
	
	\bibliographystyle{amsplain}
	\bibliography{BGS_part2.bib}

\bigskip
\noindent{\sc Keshab Chandra Bakshi} (\texttt{keshab@iitk.ac.in, bakshi209@gmail.com})\\
{\sc Satyajit Guin} (\texttt{sguin@iitk.ac.in})\\
{\sc Biplab Pal} (\texttt{biplabpal32@gmail.com, bpal21@iitk.ac.in})\\
         {\footnotesize Department of Mathematics and Statistics,\\
         Indian Institute of Technology, Kanpur,\\
         Uttar Pradesh 208016, India.}

\bigskip

\noindent{\sc Sruthymurali} (\texttt{sruthy92smk@gmail.com, sruthymurali@kannuruniv.ac.in})\\
         {\footnotesize Department of Mathematical Sciences,\\
         Kannur University,	Mangattuparamba Campus, Kannur,\\
         	Kerala 670567, India.}
         \bigskip

\end{document}